\documentclass[sn-mathphys-num]{sn-jnl}

\usepackage{graphicx}%
\usepackage{multirow}%
\usepackage{amsmath,amssymb,amsfonts}%
\usepackage{amsthm}%
\usepackage{mathrsfs}%
\usepackage[title]{appendix}%
\usepackage{xcolor}%
\usepackage{textcomp}%
\usepackage{manyfoot}%
\usepackage{booktabs}%
\usepackage{algorithm}%
\usepackage{algorithmicx}%
\usepackage{algpseudocode}%
\usepackage{listings}
\usepackage{mathtools}   
\usepackage{physics}
\usepackage{stmaryrd}    
\usepackage{tikz-cd}     
\usepackage{esint}       
\usepackage{subcaption}  
\usepackage{orcidlink}
\usepackage[most]{tcolorbox}
\usepackage{aliascnt}
\usepackage{dsfont} 

\usepackage{geometry}
\definecolor{burgundy}{rgb}{0.5, 0.0, 0.13}
\definecolor{cinereous}{rgb}{0.6, 0.51, 0.48}
\definecolor{lightsalmonpink}{rgb}{1.0, 0.6, 0.6}

\newtheoremstyle{thmstyleone}
{18pt plus2pt minus1pt}
{18pt plus2pt minus1pt}
{\itshape}
{0pt}
{\bfseries}
{}
{.5em}
{}
\theoremstyle{thmstyleone}
\numberwithin{equation}{section}
\newtheorem{theorem}{Theorem}[section]

\newaliascnt{corollary}{theorem}
\newtheorem{corollary}[corollary]{Corollary}
\aliascntresetthe{corollary}

\newaliascnt{lemma}{theorem}
\newtheorem{lemma}[lemma]{Lemma}
\aliascntresetthe{lemma}

\newaliascnt{definition}{theorem}
\newtheorem{definition}[definition]{Definition}
\aliascntresetthe{definition}

\newaliascnt{proposition}{theorem}
\newtheorem{proposition}[proposition]{Proposition}
\aliascntresetthe{proposition}

\newaliascnt{remark}{theorem}
\newtheorem{remark}[remark]{Remark}
\aliascntresetthe{remark}

\newaliascnt{notation}{theorem}
\newtheorem{notation}[notation]{Notation}
\aliascntresetthe{notation}

\newaliascnt{assumptions}{theorem}
\newtheorem{assumptions}[assumptions]{Assumptions}
\aliascntresetthe{assumptions}

\newaliascnt{example}{theorem}
\newtheorem{example}[example]{Example}
\aliascntresetthe{example}

\newcommand{\Eq}[2]{\begin{equation}\label{#1}\begin{aligned}#2 \end{aligned}\end{equation}}

\newcommand{\theo}[2]{\rbox{\begin{theorem}\label{#1} #2 \end{theorem}}}
\newcommand{\coro}[2]{\rbox{\begin{corollary}\label{#1} #2 \end{corollary}}}
\newcommand{\lem}[2]{\bbox{\begin{lemma}\label{#1} #2 \end{lemma}}}
\newcommand{\defi}[2]{\bbox{\begin{definition}\label{#1} #2 \end{definition}}}
\newcommand{\prop}[2]{\bbox{\begin{proposition}\label{#1} #2 \end{proposition}}}

\newcommand{\nota}[2]{\bbox{\begin{notation}\label{#1} #2 \end{notation}}}

\newcommand{\bbox}[1]{\begin{tcolorbox}[arc=0mm,oversize,colback=cinereous!3!white,colframe=cinereous!100!white]#1\end{tcolorbox}}
\newcommand{\rbox}[1]{\begin{tcolorbox}[arc=0mm,oversize,colback=purple!3!white,colframe=purple!100!white]#1\end{tcolorbox}}

\renewcommand{\geq}{\geqslant}
\renewcommand{\leq}{\leqslant}

\newcommand*\diag{\mathop{}\!\mathrm{diag}}

\DeclareMathOperator{\sign}{sign}

\newcommand{\eps}{\varepsilon}
\def\Xint#1{\mathchoice
   {\XXint\displaystyle\textstyle{#1}}%
   {\XXint\textstyle\scriptstyle{#1}}%
   {\XXint\scriptstyle\scriptscriptstyle{#1}}%
   {\XXint\scriptscriptstyle\scriptscriptstyle{#1}}%
   \!\int}
\def\XXint#1#2#3{{\setbox0=\hbox{$#1{#2#3}{\int}$}
     \vcenter{\hbox{$#2#3$}}\kern-.5\wd0}}

\def\dashint{\Xint-}

\title{\textbf{Correlation functions between singular values and eigenvalues}}

\author*[1,2]{\fnm{Matthias} \sur{Allard} \orcidlink{0000-0002-5682-424X}}\email{m.allard@unimelb.edu.au}

\author[1]{\fnm{Mario} \sur{Kieburg} \orcidlink{0000-0001-8864-6022}}\email{m.kieburg@unimelb.edu.au}

\affil[1]{\orgdiv{School of Mathematics and Statistics}, \orgname{University of Melbourne}, \orgaddress{\street{813 Swanston Street}, \city{Parkville, Melbourne}, \postcode{3010}, \state{Victoria}, \country{Australia}}}
\affil[2]{\orgdiv{Department of Mathematics}, \orgname{KU Leuven}, \orgaddress{\street{Celestijnenlaan 200 B bus 2400}, \city{Leuven}, \postcode{3001},  \country{Belgium}}}

\begin{document}

\abstract{Exploiting the explicit bijection between the density of singular values and the density of eigenvalues for bi-unitarily invariant complex random matrix ensembles of finite matrix size, we aim at finding the induced probability measure on $j$ eigenvalues and $k$ singular values that we coin $j,k$-point correlation measure. We  find an expression for the $1,k$-point correlation measure which simplifies drastically when assuming that the singular values follow a polynomial ensemble, yielding a closed formula in terms of the kernel corresponding to the determinantal point process of the singular value statistics. These expressions simplify even further when the singular values are drawn from a P\'{o}lya ensemble and extend known results between the eigenvalue and singular value statistics of the corresponding bi-unitarily invariant ensemble.}



\maketitle

\tableofcontents

\section{Introduction}\label{Introduction}
\subsection{State of the art}

For general complex square matrices, there exist various different decompositions. We are interested in two in particular, namely the singular value decomposition (SVD) and the Schur decomposition with which we can obtain the eigenvalues of a matrix, i.e., 

\begin{enumerate}
    \item[(i)] \textit{Singular Value Decomposition (SVD):}
    \Eq{eq:SVD}{\forall X \in \mathbb{C}^{n\times n},\ \exists\, \Sigma \in \mathbb{R}_{+,0}^n,\ U,V \in \mathrm{U}(n), \quad \textit{s.t.} \quad X=U\Sigma V}
    with $\mathrm{U}(n)$ the group of unitary matrices and $\mathbb{R}_{+,0}$ the positive real line including $0$. With the notation $\mathbb{R}_{+}$ we denote the case when we exclude $0$. The matrix $\Sigma$ is non-negative and diagonal, and its entries are the \textit{singular values} of the matrix $X$. 
    
    \item[(ii)] \textit{Schur Decomposition:}
    \Eq{eq:EVD}{\forall X \in \mathbb{C}^{n\times n},\ \exists\, z \in \mathbb{C}^n, \ t \in \mathrm{T}(n),\, U \in \mathrm{U}(n), \quad \textit{s.t.} \quad X=UztU^\dagger,}
    where $\dagger$ denotes the Hermitian conjugation and $\mathrm{T}(n)$ the group of upper unitriangular matrices. The matrix $z$ is complex and diagonal, and its entries are the \textit{eigenvalues} of the matrix $X$.

\end{enumerate}

We note that the eigenvalue decomposition in the form $X=UDU^{-1}$, with $D$ a diagonal matrix and an invertible matrix $U \in \mathrm{GL}(n,\mathbb{C})$ is not possible for every complex matrix, hence we have chosen the Schur decomposition. Every linear transformation, represented by a complex matrix $X \in G = \mathrm{GL}(n,\mathbb{C})$ can be almost entirely (up to basis transformations) characterised by its singular values $\Sigma=\mathrm{diag}(\sigma_1,\ldots,\sigma_n)$ and its eigenvalues $z=\mathrm{diag}(z_1,\ldots,z_n)$.

Both decompositions enjoy a multitude of applications, usually either the eigenvalues or the singular values alone. However, in some situations both kinds of spectra are complementary. For instance, in the study of random non-Hermitian Hamiltonian and in particular in Quantum Chaos \cite{Roccati2024,Porras2019,Herviou2019,Brunelli2023}, the physics and the chaotic behavior of the system seems to be best captured by the singular values, without making the eigenvalues redundant. The same complementarity is true in Quantum Chromodynamics~\cite{Kanazawa2011,Kanazawa2012}
 as well as in topological statistics of Hamiltonians~\cite{Braun2022,Hahn2023,Hahn2023a} or even in Time Series Analysis of time-lagged matrices~\cite{Thurner2007,Long2023,Yao2022,Loubaton2021,Bhosale2018,Nowak2017}. Born out of these motivations, we would like to address the question about the relation between the statistics of the eigenvalues and those of the singular values of a random matrix. 
 
 A few results are known, such as the Haagerup-Larson theorem~\cite{Haagerup2000} relating the limiting probability density of the eigenvalues with those of the singular values with the help of free probability techniques. A requirement of this relation has been the bi-unitary invariance of the random matrix ensemble; more details are given in the next subsection. A related result is the single ring theorem~\cite{Feinberg1997,Feinberg1997a,Guionnet2009,Rudelson2014}. Our aim is to explore more such relations for finite matrix size and higher $k$-point correlation functions.

Many standard results about the statistics of singular values and eigenvalues can be found in~\cite{Anderson2010,Akemann2015,Deift2009,Forrester2010,Byun2024}. Recent works have looked at the resulting probability density of eigenvalues of products of random matrices e.g.~\cite{Akemann2012,Akemann2013,Akemann2013a,Akemann2014,Forrester2018,Kieburg2015,Kuijlaars2014,Kuijlaars2014a,Ipsen2015,Akemann2015a,Kieburg2019}, sum of random matrices e.g.~\cite{Kuijlaars2016,Narayanan2023} and also investigated what happens to the distribution of eigenvalues when one would delete columns and rows of the matrix~\cite{Kieburg2015,Ameur2024}.

Despite this broad variety of literature on the subject, singular values and eigenvalues are seldom studied together.   From a random matrix perspective and at finite matrix size $n$, the result in~\cite{Kieburg2016} provides a bijection between the joint probability density function of the eigenvalues and the one of the singular values, under some assumptions. The related works~\cite{Kieburg2015,Kieburg2023a,Kieburg2022} bring some tools to exploit this bijection when the singular values are drawn from a particular kind of ensembles, such as polynomial ensembles~\cite{Kuijlaars2014,Kuijlaars2014a,Kieburg2019} and, more particularly, for P\'olya ensembles, which were formerly coined polynomial ensembles of derivative type~\cite{Foerster2020}. 

Let us recall that the SVD and Schur decompositions in general exhibit ambiguities.  To render the two decompositions unique, the singular values $\{\sigma_1,\ldots,\sigma_n\}$ and the moduli of the eigenvalues $\{z_1,\ldots,z_n\}$, which will be called eigenradii $\{|z_1|,\ldots,|z_n|\}$, must be ordered, and the matrices $U$ and $V$ in~\eqref{eq:SVD} and~\eqref{eq:EVD} need to be drawn from cosets. 

In general,  there exists only one equality between singular values and eigenvalues of a matrix $X\in\mathbb{C}^{n\times n}$, which is given by the modulus of the determinant
\Eq{eq:prodsev}{
|\det (X)|=|\det(z)|=\prod_{k=1}^n |z_k|=\sqrt{\det (X^\dagger X)}=\det(\Sigma)=\prod_{k=1}^n \sigma_k.
}
However, there exist various inequalities such as Weyl's inequalities~\cite{Weyl1949}. After ordering the eigenvalues  and singular values  like $|z_1| \geq |z_2| \geq \ldots \geq |z_n|$ and $\sigma_1 \geq \sigma_2 \geq \ldots \geq \sigma_n$, the simplest  Weyl inequality reads 
\Eq{eq:Weyl1}{
\prod_{k=1}^m |z_k| \leq \prod_{k=1}^m \sigma_k\quad{\rm for\ any}\ m\leq n,
}
which implies a second one
\Eq{eq:Weyl2}{
\sum_{k=1}^m |z_k| \leq \sum_{k=1}^m \sigma_k\quad{\rm for\ any}\ m\leq n.
}
Two immediate consequences follow from those two inequalities. Firstly, the largest singular value bounds the largest eigenradius from above, which is just the case $m=1$  of~\eqref{eq:Weyl1}. Secondly, the smallest  eigenradius is bounded from below by the smallest singular value as we can apply~\eqref{eq:Weyl1} for $m=1$ for the inverse matrix $X^{-1}$ if existent, otherwise there is no non-zero bound. Summarising, it is always
\Eq{eq:boundsvev}{
\sigma_1 \geq |z_1|\qquad{\rm and}\qquad \sigma_n \leq |z_n|.
}
As these relations hold for deterministic matrices, they must also hold for random matrices. Therefore, these bounds might be the source of non-trivial correlations between eigenvalues and singular values which may even survive in the limit of large matrix size.

Putting the limit of large matrix dimensions aside, we address in the present work the correlations between singular values and eigenvalues at finite matrix dimension. Thus, we aim for developing some kind of finite-dimensional counterpart of the Haagerup-Larson theorem. We postpone the question of correlations for infinite matrix dimensions to follow-up works \cite{Allard2025}.

\subsection{Main results}\label{Main results}

Assuming the probability distribution of a complex square random matrix has a density with respect to the Lebesgue measure on $\mathbb{C}^{n\times n}$, denoted by $f_G$, which does not depend on its singular vectors (right as well as left ones), then it was shown in~\cite{Kieburg2016} that the distributions of eigenvalues and singular values also have densities and there exists a linear bijection between the two densities. The property that the density of the random matrix does not depend on its singular vectors is encoded by a bi-unitary invariance of $f_G$, i.e.,
\Eq{}{
f_G(U_1 X U_2)= f_G(X)\quad{\rm for\ all}\  U_1,U_2 \in \mathrm{U}(n)\ {\rm and}\ X \in \mathrm{GL}(n,\mathbb{C}).
}
Two random matrices $X,Y \in \mathrm{GL}(n,\mathbb{C})$ are therefore equal in distribution, if they are related by $Y=U_1 X U_2$ with $U_1,U_2 \in \mathrm{U}(n)$ statistically independent of $X$ and $Y$. 
We resort to the general linear group $\mathrm{GL}(n,\mathbb{C})$ instead of $\mathbb{C}^{n\times n}$ as it is sometimes useful to guarantee the existence of an inverse $X^{-1}$. It is not problematic as $\mathrm{GL}(n,\mathbb{C})$ is dense in $\mathbb{C}^{n\times n}$ and we consider only densities so that the set of non-invertible matrices is only of measure zero. 

This impact of the bi-unitarily invariance of  random matrices should be seen in contrast to when there is no such invariance. Then, there is not much information about the relation between the two kinds of spectral statistics as the bijection between the probability distributions is lost. 

The question is then: keeping the bi-unitary invariance on $f_G$, can we find an explicit formula for the joint probability density function of the singular values \textit{and} the eigenvalues together? This turns out to be a difficult question. Especially, that the underlying probability measure will not have a density function despite that $f_G$ is a density,  due to~\eqref{eq:prodsev}.  Nonetheless, we will prove that the marginal probability measure between one eigenvalue and $k$ singular values is still a density function for a matrix size $n>1$. We will call this density the $1,k$-point correlation function whose name is reminiscent to the $k$-point correlation functions of either only eigenvalues or only singular values; see \autoref{sec:jk-densities} for a general definition. 

The derivation of explicit formulas for the $1,k$-point correlation function is one of the main goals of the present work. We will also focus on the interaction between one eigenvalue and one singular value which is captured by the $1,1$-point correlation function and in particular by the cross-covariance density function; simply defined as the difference between the $1,1$-point correlation function and the product of the respective $1$-point functions, cf. \autoref{def: cross-covariance}. Explicit formulas for all the $j,k$-point correlation functions with $j>1$ are not accessible, yet, due to some serious technical obstacles which need to be overcome. These difficulties are discussed in \autoref{Discussion}.

Due to the invertibility of $X\in\mathrm{GL}(n,\mathbb{C})$, the eigenradii and singular values are strictly positive. Actually, we will work with squared singular values and squared eigenradii to simplify the notation.

Starting from the bijective map between the probability densities of the eigenvalues and singular values for bi-unitarily invariant random matrix ensembles on $\mathrm{GL}(n,\mathbb{C})$, see~\cite{Kieburg2016}, we can derive a general expression for the $1,k$-point correlation function for $n>2$ which is summarised in the following theorem  proven in \autoref{Proof of theo:1,kpt}.

\theo{theo:1,kpt}{
Given an integer $n>2$, $k\in \llbracket 1,n\rrbracket$, contours $\mathcal{C}_j= j+i\mathbb{R}$ and $(n-1)$-dimensional vectors $\tau(j)=(1,\dots,j-1,j+1,\dots,n)$ with $j\in \llbracket 1,n\rrbracket$, where $\llbracket , \rrbracket$ denotes integer intervals. Let $f_{\rm SV}\in \mathrm{L}^1(\mathbb{R}_+^n)$ be the joint probability density of the squared singular values of a random matrix $X\in\mathrm{GL}(n,\mathbb{C})$ drawn from a bi-unitarily invariant ensemble. Then, the $1,k$-point correlation function $f_{1,k}:\mathbb{R}_+^{1+k}\to[0,\infty)$ for a squared eigenradius $r$ and $k$ pairwise distinct squared singular values $a_1,\ldots,a_k$ is
\Eq{eq:1,1pt}{
f_{1,k}(r;a_1,\ldots,a_k)=& \, \frac{1}{n} \left(\prod_{p=0}^{n-1} p! \right) \sum_{j=1}^{n}   \int_{\mathbb{R}_+^{n-k}}\left[\prod_{b=k+1}^n da_b\right]\int_{\mathcal{C}_j}\frac{ds}{2\pi i} r^{j-1-s}\\
&\times  f_{\rm SV}(a) \frac{
\det \begin{bmatrix}
    a_b^{s-1} \\
    a_b^{\tau_c(j)-1}\\
  \end{bmatrix}_{\substack{b=1,\ldots,n\\ c=1,\ldots,n-1}} }{\Delta_n(s,\tau(j))\Delta_n(a)},
}
where the determinant in the numerator should be read as follows: the first row is given by $a_b^{s-1}$ with $b=1,\ldots,n$ as the column index and the last $n-1$ rows are $a_b^{\tau_c(l)} $ with $c=1,\ldots,n-1$ as the row index. When $k=n$ no singular value is integrated out, the only remaining integral is the complex one.
Here, the $n$-dimensional Vandermonde determinant of an $n$-dimensional vector $x\in\mathbb{C}^n$ is denoted by 
\begin{equation}\label{Vandermonde}
\Delta_n(x)=\det [x_l^{p-1}  ]_{l,p=1}^n=\prod_{1\leq l<p\leq n}(x_p-x_l).
\end{equation}
The $n$-dimensional vector $(s,\tau(j))$ has $s$ as its first component and the components of $\tau(j)$ as its $(n-1)$ last entries.

The result for degenerate $a_1,\ldots,a_k$ follows from applying l'H\^opital's rule as long as $k<n$ or there is at least one pair which is distinct. For $k=n$, the limit to the fully degenerate case $a_1,\ldots,a_n\to a_1$ only exists in a weak sense.
}

The strategy to get to the $1,k$-point function is to fix $k$ squared singular values in $f_{\rm SV}$, and then use the bijection of~\cite{Kieburg2016} to get to the eigenvalues. After integrating over all eigenangles, i.e., the angles of the complex phase of the eigenvalues, and all but one eigenradius we arrive at \autoref{theo:1,kpt}.

For the case $n=1$, the induced $1,1$-point measure does not have a density, cf. \autoref{prop:rho n=1}. The case $n=2$ is given explicitly in \autoref{prop:rho n=2}, in particular~\eqref{eq:1,1-point-n2}. The proving techniques of these two results are very different than those for \autoref{theo:1,kpt} and are based on direct integration while for \autoref{theo:1,kpt} one needs to take special care of the various integrations involved.

An immediate corollary is the conditional level density of the eigenradii conditioned under fixed squared singular values $a_1,\ldots, a_n$ which can be obtained by dividing $f_{1,n}(r;a_1,\ldots,a_n)$ by $f_{\rm SV}(a)$.

\coro{cor:cond.dens}{
Under the setting of \autoref{theo:1,kpt} (especially $n>2$) and $a\in\mathbb{R}_+^n$ with pairwise distinct entries, the level density $\rho_{\rm EV}(.|a)$ of the squared eigenradii conditioned under fixed $a$ is
\begin{equation}\label{cond.dens}
\begin{split}
\rho_{\rm EV}(r|a)=&\frac{1}{n} \left(\prod_{p=0}^{n-1} p! \right) \sum_{j=1}^{n} \int_{\mathcal{C}_j}\frac{ds}{2\pi i}  r^{j-1-s} \frac{
\det \begin{bmatrix}
    a_b^{s-1} \\
    a_b^{\tau_c(j)-1}\\
  \end{bmatrix}_{\substack{b=1,\ldots,n\\ c=1,\ldots,n-1}} }{\Delta_n(s,\tau(j))\Delta_n(a)}\\
  =&\frac{1}{n}\partial_r\frac{
\det \left[\begin{array}{c|c} 0 &
    \displaystyle -\left(1-\frac{a_b}{r}\right)^{n-1}\Theta(r-a_b)  \\\hline 1 &
   \displaystyle \binom{n-1}{c-1}\left(-\frac{a_b}{r}\right)^{c-1} 
  \end{array}\right]_{b,c=1,\ldots,n} }{\det \begin{bmatrix}
    \displaystyle \binom{n-1}{c-1}\left(-\frac{a_b}{r}\right)^{c-1} \\
  \end{bmatrix}_{b,c=1,\ldots,n}}.
  \end{split}
\end{equation}
We emphasise that $b$ is the column index and $c$ the row index for both determinants and $\Theta$ is the Heaviside step function.

For degenerate squared singular values, with at least one distinct pair, one needs to apply l'H\^opital's rule. The limit to the fully degenerate case needs to be understood in the weak sense yielding the Dirac delta function
\begin{equation}
\rho_{\rm EV}(r|a_0,\ldots,a_0)=\delta(r-a_0).
\end{equation}
}

The second line can be compared with the $n=2$ result~\eqref{mu12-n2} which agree. This highlights that only the technical details of the proof are affected for this specific case but not the result. The normalisation can readily be checked, though one needs to be careful at $r=\infty$. The proper way to deal with this is to integrate over $[0,R]$ for an $R>\max_b\{a_b\}$ by expanding the first row and evaluating the Heaviside step function. Due to the derivative and the high vanishing order at $r=a_b$ the integral would be evaluated at $r=R$. After the remaining determinants are evaluated, one finds the $R$ independent normalisation $\int_0^R \rho_{\rm EV}(r|a)dr=1$ for any $R>\max_b\{a_b\}$ and, thence, for $R\to\infty$ too.

 One result one can read of from~\eqref{eq:1,1pt} is the $1$-point correlation function of the squared eigenradii $\rho_{\rm EV}$ in terms of $f_{\rm SV}$, while the one for the squared singular values $\rho_{\rm SV}$ is immediate after integrating over $r\in\mathbb{R}_+$ for $k=1$. Both are, by definition, marginal densities when integrating, respectively, over $a_1,\ldots,a_k$ or $r$  in~\eqref{eq:1,1pt}. To get a more compact expression for $\rho_{\rm EV}$, we introduce the Mellin transform  $\mathcal{M}$ on $\mathbb{R}_+$, 
 \Eq{eq:M-trans}{
\mathcal{M}f(s)=\int_0^\infty dx\ x^{s-1}f(x)
}
for an $\mathrm{L}^1(\mathbb{R}_+)$-function $f$ and $s\in\mathbb{C}$ such that the integral converges absolutely, and the spherical transform  $\mathcal{S}$ on $\mathbb{R}_+^n$,
\Eq{eq:S-trans}{
\mathcal{S}f(s)=\int_A\prod_{j=1}^n da_j\ f(a) \frac{\det [a_b^{s_c-1}  ]_{b,c=1}^n }{\Delta_n(s)\Delta_n(a)} 
}
for an $\mathrm{L}^1(\mathbb{R}_+^n)$-function $f$ and $s \in \mathbb{C}^n$ for which the integrand is Lebesgue integrable. Then, the 1-point function of the squared eigenradii $\rho_{\rm EV}$ is given by the following corollary, whose proof is a trivial identification of the terms in~\eqref{eq:1,1pt} with the definitions~\eqref{eq:M-trans} and~\eqref{eq:S-trans}.

\coro{theo:1pt}{
Consider the setting of \autoref{theo:1,kpt} apart from $n\in\mathbb{N}$ which is {\it not} necessarily larger than $2$. The 1-point correlation function of the squared eigenradii  is given by 
\Eq{eq:1point}{
\rho_{\rm EV}(r)=\frac{1}{n}\left(\prod_{p=0}^{n-1} p! \right) \sum_{j=1}^{n} r^{j-1} \mathcal{M}^{-1}\left[\mathcal{S}f_{\rm SV} (.,\tau(j))\right](r).
}
The transformation $\mathcal{M}^{-1}$ is the inverse Mellin transform~\eqref{eq:invM-trans} on $\mathbb{R}_+$, acting on the function $s\mapsto \mathcal{S}f_{\rm SV} (s,\tau(j))$, where $(s,\tau(j))$ is the same $n$-dimensional vector as in \autoref{theo:1,kpt}.
}

For $n=1$,~\eqref{eq:1point} simplifies to $\rho_{\rm EV}=\rho_{\rm SV}=f_{\rm SV}$, because the one-dimensional spherical transform reduces to the Mellin transform, i.e., $\mathcal{S}=\mathcal{M}$ for $n=1$. This is consistent with the fact the eigenradius is equal to the singular value, in this case, by~\eqref{eq:prodsev}.

Although the conditional level density~\eqref{cond.dens} of the squared radii $r$ is translucent, the results~\eqref{eq:1point} and~\eqref{eq:1,1pt} involving $f_{\rm SV}$ are not very explicit due to their generality. However, when the singular values are drawn from a \textit{polynomial ensemble}~\cite{Kuijlaars2014,Kuijlaars2014a,Kuijlaars2016,Foerster2020}, we were able to derive insightful compact formulas. The probability density of such an ensemble has the form
\Eq{eq:polynomial ensemble def}{
f_{\rm SV}(x) =\frac{ \Delta_n(x)\det\left[ w_{k-1}(x_j)  \right]_{j,k=1}^n}{n!\,\det\left[ \mathcal{M}w_{k-1}(j)  \right]_{j,k=1}^n} ,
}
where $w_0,\ldots,w_{n-1}$ are weight functions on $\mathbb{R}_+$ such that $f_{\rm SV}$ is a probability density on $\mathbb{R}_+^n$.

\begin{remark}\label{rem:bui ensemble}
    As the set of bi-unitarily invariant densities on $\mathrm{GL}(n,\mathbb{C})$ is in bijection with the set of symmetric densities on $\mathbb{R}_{+}^n$, bi-unitarily invariant ensembles are identified with the underlying ensembles of their singular values; see \eqref{eq: I_SV} and \cite[Eq.(2.22)]{Kieburg2019}. The context will therefore determine whether we refer directly to the ensemble of the singular values or the corresponding bi-unitarily invariant ensemble. 
    For instance, for the complex (Wishart-)Laguerre ensemble, the corresponding bi-unitarily invariant ensemble is the complex Ginibre ensemble. For the complex  Jacobi ensemble, it is the truncated unitary ensemble, see~\cite[Examples 3.4]{Kieburg2016}. To avoid lengthy names we stick with the names that are usually employed for the singular value statistics.

\end{remark}

 Note that an explicit expression of $\rho_{\rm EV}$ existed before for the \textit{P\'olya ensembles} (cf.~\eqref{eq:polya ensemble} and Ref.~\cite[Eq.(4.7)]{Kieburg2016}). A polynomial ensemble is a P\'olya ensemble (of multiplicative type)~\cite{Kieburg2015,Kieburg2016,Kieburg2022,Foerster2020} if there exists $w$ such that 
 \begin{equation}
 w_{k}(x)= (-x\partial_x)^{k} w(x) \in \mathrm{L}^1(\mathbb{R}_+)\qquad \forall k \in \llbracket 0,n-1\rrbracket.
 \end{equation}
To guarantee that we deal with probability measures it has been shown in~\cite{Foerster2020} that $w$ is then related to P\'olya frequency functions; thus their name. 

An advantageous property of polynomial ensembles is that their singular values follow a \textit{determinantal point process}; see~\cite{Anderson2010,Akemann2015,Deift2009,Forrester2010}. This means that the joint probability distribution $f_{\rm SV}$ can be written in the form
\Eq{kernel-polynomial}{
f_{\rm SV}(x) =\frac{1}{n!}\det\left[K(x_j,x_k)  \right]_{j,k=1}^n,
}
where $K$ is the kernel function, and all the $k$-point correlation functions have a similar form where only the size of the determinant changes. In general, $K$ is not uniquely given. Indeed, due to elementary properties of the determinant, for a non-vanishing function $g$, the kernel $[g(x_1)/g(x_2)]K(x_1,x_2)$ is also a correlation kernel for the same point process. However, we will require $K$ to be polynomial of degree $n-1$ in the second entry, which thus makes its choice unique for polynomial ensembles. Interestingly, the kernel $K$ plays also a crucial role in the correlations between the singular values and the eigenradii for polynomial ensembles, as it can be seen in \autoref{theo:poly ensemble}, proven in~\autoref{sec:proof.1.3}, and, in particular, in the structure of the cross-covariance density function, cf. \autoref{def: cross-covariance}. This theorem is understood by us as the main result of the present work as it reveals a novel and relatively simple and compact formula for the joint probability density of one squared eigenradius and $k$ squared singular values.

\theo{theo:poly ensemble}{
Let $n\in \mathbb{N}$, $n>2$, $k\in \llbracket 1,n\rrbracket$ and consider a random matrix  that is drawn from a bi-unitarily invariant ensemble on $\mathrm{GL}(n,\mathbb{C})$  having a polynomial ensemble with joint probability density~\eqref{kernel-polynomial} for the squared singular values. The $1,k$-point correlation function between one squared eigenradius and one squared singular value is given by
\begin{equation}\label{eq:1,kpt poly}
\begin{split}
&f_{1,k}(r;a_1,\ldots,a_k)\\
&=\frac{(n-k)!}{n!}\partial_r\det \left[\begin{array}{c|c}   
\displaystyle \int_0^\infty \frac{dt}{(1+t)^{n+1}}\int_0^r dv\left(1-\frac{v}{r}\right)^{n-1} K(v,-rt)  &
   \displaystyle \Hat{\Omega}(r;a_c) \\\hline \displaystyle\int_0^\infty \frac{dt}{(1+t)^{n+1}}K(a_b,-rt) & K(a_b,a_c)
  \end{array}\right]_{b,c=1}^k\\
&=\frac{(n-k)!}{(n-1)!}\int_{0}^{\infty} dt  \det\left[\begin{array}{c|c} 
    	 \displaystyle \int_0^r \frac{dv}{v} \varphi\left(\frac{v}{r},t\right) K\left(v,-rt\right)\quad & \Omega(r,a_c,t) \\ \hline
     K\left(a_b,-rt\right)\quad & K(a_b,a_c) 
\end{array}\right]_{b,c=1}^k,
\end{split}
\end{equation}
with  
\begin{equation}
\Hat{\Omega}(r;a)=\int_0^r dv\ K(v,a)\left(1-\frac{v}{r}\right)^{n-1}-\left(1-\frac{a}{r}\right)^{n-1}\Theta(r-a),
\end{equation}
 $\Theta$ the Heaviside step function,
\Eq{def:phi}{
 \varphi(x,t):= x(1-x)^{n-2}(1+t)^{-(n+2)}\biggl[\left(1-\frac{x}{n}\right)(1+t)-(1-x)\left(1+\frac{1}{n}\right) \biggl]
}
and
\Eq{def:omega}{
\Omega(r,a,t):=\int_0^r \frac{dv}{v} \varphi\left(\frac{v}{r},t\right) K(v,a) - \frac{\Theta(r-a)}{a} \varphi\left(\frac{a}{r},t\right).
}
The function $K$ is the the correlation kernel of the polynomial ensemble chosen to be a polynomial of degree $n-1$ in its second argument. 

The $1$-point functions of the squared singular values is
\Eq{}{\rho_{\rm SV}(a):=f_{0,1}(a)=\frac{1}{n}K(a,a)
}
and the one of the squared eigenradii is
\begin{equation}\label{eq:1pointreal}
\begin{split}
\rho_{\rm EV}(r):=f_{1,0}(r)=&\partial_r\int_0^\infty \frac{dt}{(1+t)^{n+1}}\int_0^r dv\left(1-\frac{v}{r}\right)^{n-1} K(v,-rt) \\
=&n \int_{0}^{\infty} dt \int_0^r \frac{dv}{v} \varphi\left(\frac{v}{r},t\right) K\left(v,-rt\right).
\end{split}
\end{equation}
}

\begin{remark}\label{rem:cov}
    The 1,k-point correlation function can be recast into the form
\Eq{}{
f_{1,k}(r; a_1,\ldots,a_k)=\rho_{\rm EV}(r)f_{0,k}(a_1,\ldots,a_k)+\mathrm{cov}_{1,k}(r;a_{1},\ldots,a_{k}),
}
where $f_{0,k}$ is, apart from the factor $(n-k)!/n!$, the standard $k$-point correlation function of the squared singular values. We coin the function $\mathrm{cov}_{1,k}(r;a_{1},\ldots,a_{k})$ as $1,k$-cross-covariance density function defined in \autoref{def: cross-covariance}, and denote $\mathrm{cov}_{1,1}$ simply by $\mathrm{cov}$. It is directly related to the covariance between an observable depending on one squared eigenradius and an observable depending on $k$ squared singular values, cf., \autoref{rem:cross-covariance}. We note that the cross-covariance density function can be seen as an analogue of the cluster functions defined in~\cite[Eq.(5)]{Dyson1962}.
\end{remark}

 For $k=1$ \autoref{theo:poly ensemble} yields 1,1-point correlation
\begin{equation}\label{eq:1,1pt poly}
\begin{split}
f_{1,1}(r; a)=&\int_{0}^{\infty} dt  \det\left[\begin{array}{c c} 
    	 \displaystyle \int_0^r \frac{dv}{v} \varphi\left(\frac{v}{r},t\right) K\left(v,-rt\right)\quad & \Omega(r,a,t) \\
    K\left(a,-rt\right) & K(a,a) \\   
\end{array}\right]=\rho_{\rm EV}(r)\rho_{\rm SV}(a)+\mathrm{cov}(r;a).
\end{split}
\end{equation}
with
\Eq{eq: cov def}{
\mathrm{cov}(r;a):=\mathrm{cov}_{1,1}(r;a)=-\int_{0}^{\infty} dt\, \Omega(r,a,t)K(a,-rt).
}

The explicit expression~\eqref{eq:1pointreal} for the $1$-point function $\rho_{\rm EV}$ of the squared eigenradii is new for a general polynomial ensemble. It is, however, known for P\'olya ensembles for which it drastically simplifies, cf.~\cite[Eq.(4.7)]{Kieburg2016}. Furthermore, the expressions in Theorem~\eqref{theo:poly ensemble} involve the function $\varphi$, which is independent of the chosen polynomial ensemble and is reminiscent of the conditional level density~\eqref{cond.dens}. 

 The expression for the $1,k$-point correlation function $f_{1,k}$ \eqref{eq:1,kpt poly}  simplifies for a P\'olya ensemble, as well. Actually, only the cross-correlation function is new for this subclass of ensembles and will be stated in the following proposition, proven in \autoref{comprho}.

\prop{coro:poly ensemble explicit}{
Let $n\in \mathbb{N}$, $n>2$.  With the same assumptions and notations as in \autoref{theo:poly ensemble}, we assume that $f_{\rm SV}$ is the joint probability density of the squared singular values of a P\'olya ensemble associated to an $n$-times differentiable weight function $w\in \mathrm{C}^n(\mathbb{R}_+)$. Then, the $1,k$-cross-covariance density function has the form
\Eq{cov.prop}{
\mathrm{cov}_{1,k}(r;a_1,\ldots,a_k)= \frac{(n-k)!}{(n-1)!}\left.\partial_\mu\det[K(a_b,a_c)+\mu\, \hat{C}(r;a_b,a_c)]_{b,c=1}^k\right|_{\mu=0}
}
with
\Eq{cov.prop2}{
\hat{C}(r;a_1,a_2)=& \sum_{\gamma=0,1} H_\gamma(r,a_2)\left[\Theta(r-a_1) \frac{1}{a_1} \Psi_\gamma\left(\frac{a_1}{r}\right)- V_\gamma(r,a_1) \right]
}
with $\Theta$ the Heaviside step function and for $\gamma=0,1$ we have employed the functions
\begin{eqnarray}
\Psi_\gamma(x)&:=& \left(\frac{1-x}{nx-1}\right)^\gamma x(1-x)^{n-2} \left(nx-1\right),\label{psi.def}\\
H_\gamma(x,y)& :=&\int_0^1 du\ q_n(yu) \partial_{u}^\gamma \left[ u^\gamma\label{H.def} \frac{\rho_{\rm EV}(xu)}{w(xu)}\right] , \\
V_\gamma(x,y) &=&  \int_0^1 du\  p_{n-1}\left(yu\right) \left(u \partial_{u}\right)^{1-\gamma} w(xu),\label{V.def}
\end{eqnarray}
where $p_{n-1}$ and $q_n$ are the bi-orthonormal pair of functions composing the kernel~\eqref{eq:polya ensemble} of  $f_{\rm SV}$ which can be expressed as
\Eq{biorthogonal}{
p_{n-1}(x) &=\sum_{c=0}^{n-1} \binom{n-1}{c} \frac{(-x)^c}{\mathcal{M}w(c+1)} \quad{\rm and}\quad q_n(x) = \frac{1}{n!}\partial_x^n[x^n w(x)]
}
according to~\cite[Lemma 4.2]{Kieburg2016}.
}

One can go even further and carry out the remaining integrals $H_\gamma(x,y)$ and $V_\gamma(x,y)$ to get a computationally efficient formula in order to create plots, cf. \autoref{lem:poly ensemble explicit}. Especially, for the classical P\'olya ensembles like Jacobi, Laguerre or Cauchy-Lorentz ensembles this is manageable. The formulation of \autoref{coro:poly ensemble explicit} might be useful for the asymptotic study $n \to \infty$, which we, however, do not address in the current work.

When $f_{\rm SV}$ is a P\'olya ensemble, the additional structure we have from a general polynomial ensemble imposes differentiability conditions on the kernel $K$ and, as a consequence, imposes continuity and differentiability conditions on the $1,k$-point correlation function $f_{1,k}$. We have analysed the analytical behaviour and proved the following conclusion in \autoref{proof coro deriv polya}.

\coro{coro:deriv Polya}{
Let $n\in \mathbb{N}$, $n>1$. With the same assumptions and notations as in \autoref{theo:poly ensemble}, if  $f_{\rm SV}$ is a P\'olya ensemble with a weight function $w\in \mathrm{C}^\infty(\sigma)$ which is smooth on the support $\sigma\subset\mathbb{R}_+$ of the $1$-point function $\rho_{\rm SV}$, then, for $n=2$, $f_{1,k}$ is discontinuous. For $n\geq 3$, it is $f_{1,k} \in \mathrm{C}^{n-3}(\sigma^{k+1})$ while it is not $(n-2)$-times continuous differentiable along  $a_j=r$ for any $j=1,\ldots,k$.
}

The present work is organized as follows. In \autoref{Preliminaries}, we present the different notations that will be used and introduce various integral transformations. Additionally, we define the $j,k$-point correlation measures and prove general expressions for the matrix sizes $n=1$ and $n=2$. We also introduce and briefly discuss polynomial and P\'olya ensembles. The proofs of the main theorems are given in \autoref{1,k-point function}. As an application and to make our results more tangible, we study the case of P\'olya ensembles, in \autoref{sec:app}. We especially give very explicit results for the Laguerre and the Jacobi ensembles (also known as Ginibre ensemble and truncated unitary ensemble, respectively). We discuss the implications of our results in \autoref{Discussion}.

\section{Preliminaries}\label{Preliminaries}
\subsection{Notations}\label{sec:notation}

For the present work, we will borrow most of the notations from~\cite{Kieburg2016}. The different matrix spaces and the corresponding measures used on them are presented in \autoref{tab:MS}. First, let us recall that given a measure on $G=\mathrm{GL}(n,\mathbb{C})$ with a density, each of the two induced measures of the singular values and of the eigenvalues have densities, too, by Tonneli's Theorem. We will denote $f_{\rm EV}: Z \to \mathbb{R}$ the density function of the eigenvalues and $f_{\rm SV}: A \to \mathbb{R}$ the density function of the squared singular values.

\begin{table}[h!] 
\begin{center}
\renewcommand{\arraystretch}{1.7}
\resizebox{\columnwidth}{!}{\begin{tabular}{ |c|c|c|} 
    \hline
    Matrix Space & Description & Reference Measure\\
    \hline \hline
    $G = \mathrm{GL}(n,\mathbb{C})$ & General linear group & $\prod_{j,k} dx_{jk}$ \\
    \hline
    $A=[\mathrm{GL}(1,\mathbb{C})/\mathrm{U}(1)]^n \cong \mathbb{R}_{+}^n$ & Group of positive definite diagonal matrices & $da=\prod_{k=1}^n da_k$ \\
    \hline 
    $Z=\mathrm{GL}(1,\mathbb{C})^n \cong \mathbb{C}^n \backslash \{0\}$ & Group of invertible complex diagonal matrices & $dz=\prod_{k=1}^n dz_k$ \\
    \hline 
    $\mathrm{U}(n)$ & Group of unitary matrices & $d\mu_H(u)=$ normalized Haar measure \\
    \hline 
    $\mathrm{T}(n)$ & Group of upper unitriangular matrices & $dt=\prod_{j>k} dt_{jk}$ \\
    \hline    
\end{tabular}}
\end{center}
\caption{\footnotesize\textsc{Matrix Spaces and Reference Measures} (notation adapted from \protect\cite[Table 1]{Kieburg2016}). Here, $dx$ denotes the Lebesgue measure on $\mathbb{R}$ if $x$ is a real variable and the Lebesgue measure $dx = d\Re{x} d \Im{x}$ on $\mathbb{C}$ if $x$ is a complex variable.\label{tab:MS}}
\end{table}

By abuse of notation, we will identify vectors of eigenvalues, squared eigenradii and squared singular values with diagonal matrices out of convenience. The squared singular values  $a=\diag(a_1,\ldots,a_n)$ and the eigenvalues $z=\diag(z_1,\ldots,z_n)$ will be unordered. 

\subsection{Harmonic Analysis}

Our methods are based on harmonic analysis tools and the bijection proven in~\cite[Theorem 3.1]{Kieburg2016}, see \autoref{theo:R-transform}. Thus, we will briefly recall the corresponding transforms and introduce our notation for those.

We start with the \textit{Mellin transform} for a measurable function $f$ on $\mathbb{R}_+$,  which is defined in~\eqref{eq:M-trans}. When $f$ is a probability density, the normalisation is given by $\mathcal{M}f(1)=1$. The Mellin transform is only defined for those $s \in \mathbb{C}$ such that the integral exists (in the Lebesgue sense). In particular, if $f \in \mathrm{L}^1(\mathbb{R}_+)$, the Mellin transform is defined at least on the line $\mathcal{C}_1:=1+i\mathbb{R}$. Let $\mathrm{L}^1(\mathbb{R}_+)$ be the space of Lebesgue integrable functions on $\mathbb{R}_+$.
By the Mellin inversion theorem, e.g., see~\cite[Lemma 2.6]{Kieburg2016}, $\mathcal{M}:\mathrm{L}^1(\mathbb{R}_+)\to \mathcal{M}\mathrm{L}^1(\mathbb{R}_+)$ is bijective and the Mellin inversion formula can be  given by the limit
\Eq{eq:invM-trans}{
\mathcal{M}^{-1}[\mathcal{M}f](x):=\lim_{\varepsilon \to 0} \int_{\mathcal{C}_1}\frac{ds}{2\pi i} \zeta(\varepsilon \Im{s}) x^{-s}\mathcal{M}f(s)=f(x),
}
with the regularisation $\zeta$ defined as in~\cite[Eq.(2.40)]{Kieburg2016}, where it is denoted $\zeta_1$; in particular it is
\Eq{eq:reg fct}{
 \zeta(s):=\frac{\cos(s)}{1-4 s^2/\pi^2}.
}
The function $\zeta$ guarantees the absolute integrability and makes the Mellin transformation bijective. It can be dropped when $\mathcal{M}f$ is absolutely integrable on $\mathcal{C}_1$.

 We also need the multivariate version of the Mellin transform, which can be defined using the tensor product $\otimes$,
\Eq{eq:multiM-trans}{
\mathcal{M}^{\otimes n}f(s)=\int_A da\, f(a) \prod_{k=1}^n a_k^{s_k-1}.
}
As we are working with  densities symmetric under permutation of their arguments, we need transformations that preserve the symmetry. Particularly, we assume $f\in\mathrm{L}^{1,\rm SV}(A)$ where $\mathrm{L}^{1,\rm SV}(A)$ is the space of symmetric Lebesgue integrable functions on $A$ in which also the joint probability densities for the squared singular values can be found, thus, the chosen notation. Therefore, we can go over to the symmetrised version of the multivariate Mellin transform given by
\Eq{eq:SymmultiM-trans}{
\mathcal{M}_{\mathrm{S}}f(s)=\frac{1}{n!} \sum_{\sigma \in \mathrm{S}_n} \mathcal{M}^{\otimes n}f(\sigma(s))=\frac{1}{n!} \int_A da\,f(a)\, \mathrm{Perm} [a_j^{s_k-1}]_{j,k=1}^n  ,
}
with $\mathrm{S}_n$ the finite symmetric group of permutations of $n$ elements,  $\sigma(s)=(\sigma(s_1),\ldots,\sigma(s_n))$ and the permanent
\Eq{}{
\mathrm{Perm}[x_{jk}]_{j,k=1}^n &=\sum_{\sigma \in \mathrm{S}_n}\prod_{j=1}^n x_{j\sigma(j)}.
}
The symmetrized inverse Mellin transform~\cite[Eq.(2.39)]{Kieburg2016} is, then, given by
\Eq{inv.Mellin}{
\mathcal{M}_{\mathrm{S}}^{-1}[\mathcal{M}_{\mathrm{S}}f](x)=\frac{1}{n!} \lim_{\varepsilon \to 0} \int_{\mathcal{C}(n)} \left[\prod_{k=1}^n \frac{ds_k}{2\pi i}\zeta(\varepsilon \Im{s_k})\right] \mathrm{Perm} [x_j^{-s_k}]_{j,k=1}^n \mathcal{M}_{\mathrm{S}}f(s) ,
}
with $\mathcal{C}(n)=\bigtimes_{k=1}^n \mathcal{C}_k$, the Cartesian product of elementary contours $\mathcal{C}_k= k+i\mathbb{R}$ which are straight lines parallel to the imaginary axis going from $k-i\infty$ to $k+i\infty$.

Another important multivariate integral transformation is the \textit{spherical transform} $\mathcal{S}: \mathrm{L}^{1,\rm SV}(A) \to \mathcal{S}\mathrm{L}^{1,\rm SV}(A)$ defined in~\eqref{eq:S-trans}. We use the notation $\mathcal{S}\mathrm{L}^{1,\rm SV}(A)$ to emphasize that it is the image space of $\mathcal{S}$ with respect to the domain $\mathrm{L}^{1,\rm SV}(A)$. Note that $\mathcal{S}$ preserves the permutation symmetry of $f$ in its arguments.

These Mellin and spherical transforms are the building blocks for the singular value--eigenvalue transformation, coined \textit{SEV transform} $\mathcal{R}$, defined in~\cite[Theorem 3.1]{Kieburg2016}. It is a bijective map between the set of symmetric densities on the squared singular values, $\mathrm{L}^{1,\rm SV}(A)$, and  the set of induced densities of eigenvalues of bi-unitarily invariant matrix ensembles denoted by $\mathrm{L}^{1,\rm EV}(Z)$. We recall the theorem here for convenience.

\theo{theo:R-transform}{
\textup{(\cite[Theorem 3.1]{Kieburg2016})} Let $\mathcal{C}(n)$ be the $n$-dimensional contour in~\eqref{inv.Mellin}. The map $\mathcal{R}: \mathrm{L}^{1,\rm SV}(A)\to \mathrm{L}^{1,\rm EV}(Z)$ from the joint densities of the squared singular values to the
joint densities of the eigenvalues induced by the bi-unitarily
invariant signed densities is bijective and has the explicit integral representation 
\Eq{eq:SEV}{
       f_{\rm EV}(z) &=  \mathcal{R}f_{\rm SV}(z) =\frac{\prod_{j=0}^{n-1} j!}{(n!)^2 \pi^n} |\Delta_n(z)|^2 \mathcal{M}_{\mathrm{S}}^{-1}\mathcal{S}f_{\rm SV}(|z|^2)\\
        &= \frac{\prod_{j=0}^{n-1} j!}{(n!)^2 \pi^n} |\Delta_n(z)|^2\lim_{\varepsilon \to 0} \int_{\mathcal{C}(n)}\left[\prod_{k=1}^n \frac{ds_k}{2\pi i}\zeta(\varepsilon \Im{s_k})\right] \mathrm{Perm}[|z_b|^{-2 s_c}  ]_{b,c=1}^n\\
 & \times \int_A \left[ \prod_{j=1}^{n} \frac{da_j}{a_j} \right]\ f_{\rm SV}(a) \frac{\det [a_b^{s_c}  ]_{b,c=1}^n }{\Delta_n(s)\Delta_n(a)}  ,
         }
         where $|z|^2=\diag(|z_1|^2,\ldots,|z_n|^2)$. Especially, $\mathrm{L}^{1,\rm EV}(Z)=\mathcal{R}\mathrm{L}^{1,\rm SV}(A)$.
}

The explicit integral representation of the inverse map $\mathcal{R}^{-1}$ can be found in~[Eq.(3.4)]\cite{Kieburg2016}. Let us underline that the eigenangles only appear in the factor $|\Delta_n(z)|^2$. The bi-unitary invariance of the random matrix $X \in G$ implies that its spectrum is isotropic which, in turn, implies that the arithmetic mean of the eigenangles should be uniformly distributed on the interval $[0,2\pi]$. The differences of the eigenangles are, however, not uniformly distributed.

The linear integral transformations linking the different function spaces involved in \autoref{theo:R-transform} can be represented in the following commutative diagram which is a reduced version of~\cite[Eq.(3.1)]{Kieburg2016},
\begin{equation}\label{diagram}
\begin{tikzcd}[row sep=4em,column sep=4em]
\mathrm{L}^{1,\rm BU}(G) \arrow{r}{\mathcal{I}_{\rm SV}} \arrow[swap]{d}{\mathcal{I}_{\rm EV}} & \mathrm{L}^{1,\rm SV}(A) \arrow{ld}{ \mathcal{R}} \arrow{d}{\mathcal{S}} \\
\mathrm{L}^{1,\rm EV}(Z)   & \mathcal{S}\mathrm{L}^{1,\rm SV}(A) \arrow{l}{\mathcal{Z}}
\end{tikzcd}
\end{equation} 
where
\Eq{}{
\mathrm{L}^{1,\rm BU}(G):=\{ \,f_G \in \mathrm{L}^{1}(G)\, |\, f_G \,\, \text{is bi-unitarily invariant  on}\,\,  G \,\}.
}
The transformation $\mathcal{Z}: \mathcal{S}\mathrm{L}^{1,\rm SV}(A) \to \mathrm{L}^{1,\rm EV}(Z)$ is defined as
\Eq{}{
\mathcal{Z}\widehat{f}(z)=\frac{\prod_{j=0}^{n-1} j!}{(n!)^2 \pi^n} |\Delta_n(z)|^2 \mathcal{M}_{\mathrm{S}}^{-1}\widehat{f}(z)=f_{\rm EV}(z),
}
while $\mathcal{I}_{\rm SV}:\mathrm{L}^{1,\rm BU}(G) \to \mathrm{L}^{1,\rm SV}(A) $ consists of a singular value decomposition~\eqref{eq:SVD} and integrating out the unitary matrices $U$ and $V$ with respect to the corresponding Haar measure. Explicitly, it is
\Eq{eq: I_SV}{
\mathcal{I}_{\rm SV}f_G(a)=\left(\frac{\pi^{n^2}}{n!\,\left(\prod_{k=0}^{n-1}k! \right)^2}\right)\Delta_n(a)^2f_G(\sqrt{a})=f_{\rm SV}(a)
}
and one can thus identify the bi-unitarily invariant ensemble with the corresponding ensemble of its singular values via the relation $f_G=\mathcal{I}_{\rm SV}^{-1}f_{\rm SV}$; cf.,  \cite[Eq.(2.22)]{Kieburg2019}.
The transformation $\mathcal{I}_{\rm EV}:\mathrm{L}^{1,\rm BU}(G) \to \mathrm{L}^{1,\rm EV}(Z)$ consists of the Schur decomposition~\eqref{eq:EVD} and integrating out the Haar distributed unitary matrix $U$ as well as the upper triangular matrix, i.e.,
\Eq{eq: I_EV}{
\mathcal{I}_{\rm EV}f_G(z) &=\left(\frac{1}{n!}\prod_{k=0}^{n-1}\frac{\pi^{k}}{k!} \right)|\Delta_n(z)|^2 \left(\prod_{k=1}^{n} |z_k|^{2(n-k)} \right)\int_{\mathrm{T}(n)}dt\ f_G(zt) = f_{\rm EV}(z).
}
This latter transformation is surprisingly invertible, despite the integral over $t$, as shown in~\cite{Kieburg2016}. 

It is worthwhile to stress that any symmetric probability density function on $A$ can be traced back uniquely to a probability density of a given bi-unitarily invariant ensemble on $G$. The simplest way is  to build a corresponding bi-unitarily invariant matrix multiplying the matrix $a=\diag(a_1,\ldots,a_n)$ on the right and on the left by two independent Haar distributed unitary matrices. Unfortunately, not every  symmetric probability density function on $Z$, can be seen as the marginal distribution of bi-unitarily invariant random matrix ensemble after employing a Schur decomposition. Applying $\mathcal{R}^{-1}$ on an arbitrary symmetric probability density on $Z$ can give a signed density on $A$. This is why $\mathrm{L}^{1,\rm EV}(Z)$ is strictly a subset of all symmetric densities on $Z$ when $n>1$ and, hence, why we identify the bi-unitarily invariant ensemble with the ensemble of its singular values and not the one of its eigenvalues.

Our goal is to exploit \autoref{theo:R-transform} to explore the relationship between squared singular values and squared eigenradii. We would like to find the joint probability measure on both the squared eigenradii and squared singular values along with the induced marginal measures. Those are not necessarily densities as, in the case of the joint measure, equation~\eqref{eq:prodsev} imposes a strong constraint which might give the induced measure a component of a Dirac delta measure.

\subsection{Correlation functions}\label{sec:jk-densities}

We denote the expected value of a measurable function $\phi:G\to \mathbb{C}$ on $G$ by 
\Eq{eq:expected value}{
\mathbb{E}[\phi(X)]:=\int_G dX\ \phi(X)f_G(X).
}
With the help of this notation we define the $j,k$-point correlation measures in a weak topological sense.
\defi{def:j,k pt meas}{ 
\textup{($j,k$-point correlation measure/function)} Let $n \in \mathbb{N}$, $j,k\in \llbracket 0,n \rrbracket$. Let $f_G$ be the probability density function of the random matrix $X\in G=\mathrm{GL}(n,\mathbb{C})$. Denoting the squared eigenradii of $X$ by $\{r_l(X) \}_{l=1}^n$ and its squared singular values by $\{a_l(X) \}_{l=1}^n$, then the $j,k$-point correlation measure $\mu_{j,k}$ is defined weakly by the relation 
\Eq{eq:j,k-point pdf}{
&\mathbb{E}\bigg[\frac{1}{(n!)^2} \sum_{(\sigma,\pi)\in \mathrm{S}_n^2} \phi_{\rm EV}(r_{\sigma(1)}(X),\ldots,r_{\sigma(j)}(X))\phi_{\rm SV}( a_{\pi(1)}(X),\ldots,a_{\pi(k)}(X) ) \bigg]\\
=&\,\int_{\mathbb{R}_+^{j+k}}d\mu_{j,k}(r_1,\ldots,r_j;a_1,\ldots,a_k)\phi_{\rm EV}(r_1,\ldots,r_j)\phi_{\rm SV}(a_1,\ldots,a_k)  
}
for any Schwartz function $\phi_{\rm EV}\in\mathcal{S}\left(\mathbb{R}_+^{j}\right)$ and continuous bounded function $\phi_{\rm SV} \in \mathrm{C}_b\left(\mathbb{R}_+^{k}\right)$, where $\mathrm{S}_n$ is the symmetric group permuting $n$ elements. The induced probability measures $\mu_{j,k}$, on $j$ squared eigenradii and $k$ squared singular values of the random matrix $X$ are called $j,k$-point correlation measures. 
If the $j,k$-point correlation measure $\mu_{j,k}$ has a density with respect to the Lebesgue measure, the density will be denoted $f_{j,k}$ and will be called the $j,k$-point correlation  function.
}

We note that $\mu_{j,k}$ is indeed a probability measure due to the normalization factor $1/(n!)^2$; see \autoref{normalisation jk pt fct}.

The choice of different test functions for the squared singular values and the squared eigenradii is motivated solely by technical considerations, as, for the purpose of the proofs, additional analytical properties are needed for the function on the eigenradii. We believe, however, that this distinction is of minor significance, since both bounded continuous functions and Schwartz functions are dense in all kinds of function spaces. In particular, any test function that does not factor into separate functions of the squared singular values and the squared eigenradii can be understood as a weak limit---or another appropriate type of limit---of linear combinations of product functions of the form $\phi_{\rm EV}\phi_{\rm SV}$.

\begin{remark}
 If $j=0$ or $k=0$, we get the marginal probability measure of only $k$ squared singular values or $j$ squared eigenradii, respectively. By definition we set $\mu_{0,0}=1$ so that it is consistent with $\int_{\mathbb{R}_+} d\mu_{0,1}(a)=\int_{\mathbb{R}_+} d\mu_{1,0}(r)=\mu_{0,0}=1$.
\end{remark}

The following Lemma is rather helpful in relating the definition above with the SEV transform~\eqref{eq:SEV}.

\lem{}{
Let $z(X)$ be  the diagonal matrix of eigenvalues  and $a(X)$ comprises the squared singular values of $X\in G$. Additionally, let $f:Z\times A\to\mathbb{C}$ such that $g\in\mathrm{L}^{1}(G)$ with $g(X)=f(z(X),a(X))$. Then, we have
\begin{equation}
\begin{split}\label{relation.eig.sing}
\int_G dX\ f(z(X),a(X))=& \frac{\pi^{n(n-1)}}{(n!)^3\prod_{j=0}^{n-1} j!} \int_Z dz |\Delta_n(z)|^2\lim_{\varepsilon \to 0} \int_{\mathcal{C}(n)}\left[\prod_{k=1}^n \frac{ds_k}{2\pi i}\zeta(\varepsilon \Im{s_k})\right] \\
 & \times \mathrm{Perm}[|z_b|^{-2 s_c}  ]_{b,c=1}^n\int_A da\, f(z,a) \frac{\Delta_n(a)\det [a_b^{s_c-1}  ]_{b,c=1}^n }{\Delta_n(s)}.
 \end{split}
\end{equation}
}

We underline that $f$ is not the joint probability density of the eigenvalues and the squared singular values but some general integrable function.

\begin{proof}
The main idea is to first decouple the integral over the first $n$ arguments of $f$ from the integral over the triangular matrix which appears when performing a Schur decomposition~\eqref{eq:EVD}. Then, we can make use of the commutative diagram~\eqref{diagram}, essentially only of the triangle with the corners $\mathrm{L}^{1,\rm BU}(G)$, $\mathrm{L}^{1,\rm SV}(A)$ and $\mathrm{L}^{1,\rm EV}(Z)$. The advantage is that the complex eigenvalues $z$ are fixed in this part of the diagram.

In the first step, we perform the Schur decomposition~\eqref{eq:EVD}, i.e., $X=UztU^\dagger$. As the squared singular values are bi-unitarily invariant functions one can perform the integration over the unitary group and gets
\begin{equation}
\begin{split}\label{change.bi}
\int_G dX\ f(z(X),a(X))=& \int_Z dz\ \left(\frac{1}{n!}\prod_{k=0}^{n-1}\frac{\pi^{k}}{k!} \right)|\Delta_n(z)|^2 \left(\prod_{k=1}^{n} |z_k|^{2(n-k)} \right)\\
&\times\int_{\mathrm{T}(n)}dt f(z,a(zt)),
 \end{split}
\end{equation}
cf. Eq.~\eqref{eq: I_EV}. When considering the integral over $t$ with fixed $z$, we notice that this is the operator $\mathcal{I}_{\rm EV}$, which is, on the other hand, equal to $\mathcal{R}\circ\mathcal{I}_{\rm SV}$, see the commutative diagram~\eqref{diagram}. We underline that the SEV transform $\mathcal{R}$ also applies for general Lebesgue integrable functions and not only probability densities due to its linear nature as an operator. As required, we assumed that $g\in\mathrm{L}^{1}(G)$ with $g(X)=f(z(X),a(X))$ implying that $f$ is Lebesgue integrable in the last $n$ entries for almost all $z\in Z$ with respect to the reference measure  $|\Delta_n(z)|^2dz$. Plugging in~\eqref{eq: I_SV} and~\eqref{eq:SEV} we arrive at the assertion.
\end{proof}

\begin{remark}\label{rem:identification}
Considering \autoref{def:j,k pt meas} and the joint probability density of the squared singular values $f_{\rm SV}\in\mathrm{L}^{1,\rm SV}(A)$, especially Eq.~\eqref{eq: I_SV}, we can choose, in particular,
\begin{equation}\label{ident.test}
\begin{split}
f(z,a)=&\frac{\left(\prod_{l=0}^{n-1}l! \right)^2}{\pi^{n^2}n!}\frac{f_{\rm SV}(a)}{\Delta_n(a)^2}\sum_{(\sigma,\pi)\in \mathrm{S}_n^2}\phi_{\rm EV}\left(|z_{\sigma(1)}|^2,\ldots,|z_{\sigma(j)}|^2\right)\phi_{\rm SV}\left( a_{\pi(1)},\ldots,a_{\pi(k)} \right) .
 \end{split}
\end{equation}
\end{remark}

In the case $j=k=1$ with $n>1$, it will be shown that $\mu_{1,1}$ admits a density $f_{1,1}\in \mathrm{L}^1(\mathbb{R}^2_+)$, that will therefore be called the $1,1$-point correlation function between one squared eigenradius and one squared singular value. 

To illustrate the definition of $j,k$-point correlation measures, we consider the simplest cases of $n=1,2$ where we concentrate only on mixed correlation measures, meaning $j,k>0$. For this purpose let us introduce the Dirac distribution $\delta$, see~\cite{Schwartz1966}, which acts on any continuous function $\phi \in \mathrm{C}^0(\mathbb{R})$ as
\Eq{eq:dirac}{
\int_{\mathbb{R}}\phi(x)\delta(x-x_0)dx=\phi(x_0) \quad \textrm{for almost all}\,\, x_0 \in \mathbb{R}.
}
Then, we have the following trivial result for $n=1$.

\prop{prop:rho n=1}{
\textup{(The case $n=1$)} Let $f_G\in\mathrm{L}^{1,\rm BU}(G)$ be a probability density for $n=1$. Then, the induced probability measure $\mu_{1,1}$ on the squared singular value and the squared eigenradius is given by
\Eq{}{
d\mu_{1,1}(r_1;a_1)=\pi f_G(\sqrt{r_1})\delta(r_1-a_1)dr_1 da_1.}
}

\begin{proof}
Taking $\phi \in \mathrm{C}_b\left(\mathbb{R}_+^2\right)$, 
\Eq{}{
\mathbb{E}\bigg[  \phi(r_1(X); a_1(X)) \bigg]=\int_G dX\ \phi(r_1(X); a_1(X))f_G(X).
}
We proceed with a Schur decomposition~\eqref{eq:EVD}, which is the change to polar coordinates $X=\sqrt{r_1}e^{i\theta_1}$ for $n=1$. The measure becomes $dX=\sqrt{r_1}d\sqrt{r_1}d\theta_1=\frac{1}{2}dr_1 d\theta_1$. We use the fact that $r_1=a_1$ by~\eqref{eq:prodsev}, and integrate out $\theta_1$ which is uniformly distributed by the bi-unitary invariance of $X$. One then gets,
\Eq{}{
\mathbb{E}\bigg[  \phi(r_1(X); a_1(X)) \bigg]=\pi\int_{\mathbb{R}_+}dr_1\ \phi(r_1; r_1)f_G(\sqrt{r_1}).
}
On the other hand,~\eqref{eq:j,k-point pdf} for $\mu_{1,1}$ reads
\Eq{}{
\mathbb{E}\bigg[  \phi(r_1(X); a_1(X)) \bigg]=\int_{\mathbb{R}_+^2}d\mu_{1,1}(r_1;a_1)\ \phi(r_1; a_1).
}
Identification yields the claim.
\end{proof}

The case $n=2$ is richer with $j,k$-point correlation measures as we have now four measures with mixed statistics in the squared eigenradii and squared singular values compared to a single one for $n=1$.

\prop{prop:rho n=2}{
\textup{(The case $n=2$)} Let $f_{\rm SV}\in\mathrm{L}^{1,\rm SV}(A)$ be the joint probability density of the squared singular values for $n=2$.
Then, the $2,2$-point correlation measure is
\Eq{mu22-n2}{
d\mu_{2,2}(r_1,r_2;a_1,a_2)=&\Theta\left(\max\{a_1,a_2\}-\max\{r_1,r_2\}\right)\Theta\left(\min\{r_1,r_2\}-\min\{a_1,a_2\} \right)\\
&\times \frac{f_{\rm SV}(a_1, a_2)}{2|a_1-a_2|}\left(r_1+r_2\right)\, \delta(r_1 r_2-a_1 a_2) dr_1dr_2 da_1 da_2,
}
where $\Theta$ is the Heaviside step function.
The $2,1$-point correlation function is given by
\Eq{mu21-n2}{
f_{2,1}(r_1,r_2;a_1)=\left[\Theta\left(a_1-\max\{r_1,r_2\}\right)+\Theta\left(\min\{r_1,r_2\}-a_1 \right)\right]\frac{f_{\rm SV}\left( a_1, r_1 r_2/a_1 \right)}{2|a_1^2-r_1 r_2|}(r_1+ r_2),
}
for almost all $r_1,r_2,a_1>0$ and the $1,2$-point correlation  function by
\Eq{mu12-n2}{
f_{1,2}(r_1;a_1,a_2)=\Theta\left(\max\{a_1,a_2\}-r_1\right)\Theta\left(r_1-\min\{a_1,a_2\} \right)\frac{f_{\rm SV}(a_1, a_2)}{2|a_1-a_2|}\left(1+ \frac{a_1 a_2}{r_1^2}\right)
}
for almost all $r_1,a_1,a_2>0$.
The $1,1$-point correlation function is then
\begin{equation}\label{eq:1,1-point-n2}
\begin{split}
f_{1,1}(r_1;a_1)=&\int_{0}^{\infty} da_2f_{1,2}(r_1;a_1,a_2)=\int_{0}^{\infty} dr_2 f_{2,1}(r_1,r_2;a_1)\\
=&\left[\Theta(a_1-r_1)\int_0^{r_1}+\Theta(r_1-a_1)\int_{r_1}^\infty\right]da_2 \frac{f_{\rm SV}(a_1, a_2)}{2|a_1-a_2|}\left(1+ \frac{a_1 a_2}{r_1^2}\right).
\end{split}
\end{equation}

}

Unfortunately, the marginal density $f_{1,1}$ cannot be simplified much further unless one resorts to subclasses of ensembles. For instance, polynomial ensembles have the joint probability density of the squared singular values
\begin{equation}
f_{\rm SV}(a_1, a_2)=\frac{(a_2-a_1)[w_0(a_1)w_1(a_2)-w_1(a_1)w_0(a_2)]}{2[\mathcal{M}w_0(1)\mathcal{M}w_1(2)-\mathcal{M}w_0(2)\mathcal{M}w_1(1)]},
\end{equation}
which implies the $1,1$-point correlation function
\begin{equation}\label{Polyn.n=2}
\begin{split}
&f_{1,1}(r_1;a_1)\\
=&\frac{\Theta\left(r_1-a_1 \right)}{4}\int_{r_1}^\infty da_2\frac{w_0(a_1)w_1(a_2)-w_1(a_1)w_0(a_2)}{\mathcal{M}w_0(1)\mathcal{M}w_1(2)-\mathcal{M}w_0(2)\mathcal{M}w_1(1)}\left(1+ \frac{a_1 a_2}{r_1^2}\right)\\
&-\frac{\Theta\left(a_1-r_1\right)}{4}\int_0^{r_1}\frac{w_0(a_1)w_1(a_2)-w_1(a_1)w_0(a_2)}{\mathcal{M}w_0(1)\mathcal{M}w_1(2)-\mathcal{M}w_0(2)\mathcal{M}w_1(1)}\left(1+ \frac{a_1 a_2}{r_1^2}\right)\\
=&\frac{\Theta\left(r_1-a_1 \right)}{4}\frac{w_0(a_1)[\mathcal{M}w_1(1)+a_1\mathcal{M}w_1(2)/r_1^2]-w_1(a_1)[\mathcal{M}w_0(1)+a_1\mathcal{M}w_0(2)/r_1^2]}{\mathcal{M}w_0(1)\mathcal{M}w_1(2)-\mathcal{M}w_0(2)\mathcal{M}w_1(1)}\\
&-\frac{1}{4}\frac{w_0(a_1)[\tilde{w}_{1,r_1}(1)+a_1\tilde{w}_{1,r_1}(2)/r_1^2]-w_1(a_1)[\tilde{w}_{0,r_1}(1)+a_1\tilde{w}_{0,r_1}(2)/r_1^2]}{\mathcal{M}w_0(1)\mathcal{M}w_1(2)-\mathcal{M}w_0(2)\mathcal{M}w_1(1)},
\end{split}
\end{equation}
where $\tilde{w}_{j,r_1}$ is the incomplete Mellin transform of $w_{j}$, see~\eqref{incMel.def} below.

\begin{proof}
We start from~\eqref{change.bi} with the choice~\eqref{ident.test}. Without loss of generality we can assume that the test function $\phi$ is symmetric in its first two arguments as well as its last two ones so that the sum becomes trivial and yields a factor of $4$ cancelling with the combinatorial factor in front of the sum. For $n=2$, we use the following relation
\begin{equation}
a_1+a_2=\tr (zt(zt)^\dagger)=r_1(1+|\tau|^2)+r_2,
\end{equation}
where $\tau$ is the complex number in the off-diagonal of $t$. Plugging in $a_2=r_1r_2/a_1$ originating from the identity~\eqref{eq:prodsev}, we have
\begin{equation}
|\tau|^2=\frac{a_1}{r_1}+\frac{r_2}{a_1}-1-\frac{r_2}{r_1}
\end{equation}
which has either two or no solutions in $a_1$. The situation of no solution corresponds to $a_1\in(\min\{r_1,r_2\},\max\{r_1,r_2\})$ as, then, the right hand side is negative while the left hand side is non-negative. The case $a_1\leq\min\{r_1,r_2\}$ corresponds to the solution with the ordering $a_1<a_2$ while $a_1\geq\max\{r_1,r_2\}$ relates to $a_1>a_2$. Both branches map to the very same $|\tau|^2\geq0$ so that the substitution is not bijective. Since the situation must be invariant under swapping $a_1$ and $a_2$ the two contributions yield the very same weight meaning it yields a factor of $1/2$.

Returning to~\eqref{change.bi}, we perform a polar decomposition of $\tau$ and substitute $|\tau|$ by $a_1$. Afterwards, we integrate over the complex phases of $\tau$, $z_1$ and $z_2$. Then, we arrive at
\begin{equation}
\begin{split}
&\mathbb{E}\bigg[\phi(r_{1}(X),r_{2}(X); a_{1}(X),a_{2}(X) ) \bigg]\\
=&\frac{\pi^4}{4} \int_{\mathbb{R}_+^2} dr_1dr_2 (r_1+r_2)r_1\left[\int_0^{\min\{r_1,r_2\}}+\int_{\max\{r_1,r_2\}}^\infty\right]da_1\left|\frac{1}{r_1}-\frac{r_2}{a_1^2}\right|\\
&\times  \frac{2}{\pi^{4}}\frac{f_{\rm SV}\left(a_{1},r_1r_2/a_1\right)}{(a_1-r_1r_2/a_1)^2}\phi\left(r_1,r_2; a_{1},\frac{r_1r_2}{a_1}\right),
 \end{split}
\end{equation}
where we used Eq.~\eqref{ident.test} when plugging in the considered setting. From this equation and \autoref{def:j,k pt meas} we can read off
\Eq{}{
&d\mu_{2,2}(r_1,r_2;a_1,a_2)\\
=&\biggl[\frac{r_1+r_2}{a_1}\, \delta\left(\frac{r_1 r_2}{a_1}- a_2\right)\left[\Theta\left(a_1-\max\{r_1,r_2\}\right)+\Theta\left(\min\{r_1,r_2\}-a_1 \right)\right]\\
&+\frac{r_1+r_2}{a_2}\, \delta\left(\frac{r_1 r_2}{a_2}- a_1\right)\left[\Theta\left(a_2-\max\{r_1,r_2\}\right)+\Theta\left(\min\{r_1,r_2\}-a_2 \right)\right]\biggl]\\
&\times\frac{f_{\rm SV}(a_1, a_2)}{4|a_1-a_2|} dr_1dr_2 da_1 da_2,
}
where we have symmetrised in $a_1$ and $a_2$ as we consider unordered squared singular values.
After applying the standard rules of the Dirac delta distribution and the Heaviside step function we find~\eqref{mu22-n2}.

Claims~\eqref{mu21-n2},~\eqref{mu12-n2}, and~\eqref{eq:1,1-point-n2} can be readily obtained by setting $\phi(r_1,r_2;a_1,a_2)=[\psi(r_1,r_2;a_1)+\psi(r_1,r_2;a_2)]/2$, $\phi(r_1,r_2;a_1,a_2)=[\psi(r_1;a_1,a_2)+\psi(r_2;a_1,a_2)]/2$ or $\phi(r_1,r_2;a_1,a_2)=[\psi(r_1;a_1)+\psi(r_2;a_1)+\psi(r_1;a_2)+\psi(r_2;a_2)]/4$, respectively, with $\psi$ some other test function.
\end{proof}

We would like to turn to the $k$-point correlation function. When one is interested in marginal densities, it is suitable to consider
\Eq{eq:k-point function}{
f_k(x_1,\ldots,x_k):=f_{0,k}(x_1,\ldots,x_k)=\int_{\mathbb{R}_+^{n-k}}f(x_1,\ldots,x_n)dx_{k+1}\ldots dx_n
}
with $k\in \llbracket 1,n \rrbracket$ and $f\in \mathrm{L}^{1}(A)$.

\begin{remark}\label{normalisation jk pt fct}
This differs from the $k$-point correlation function $R_k$, used to study determinantal point processes, by a combinatorial factor~\cite{Akemann2015,Forrester2010} reminiscent of the permutation symmetry of the density in its arguments, 
\Eq{eq: R_k}{
R_k(x_1,\ldots,x_k):=\frac{n!}{(n-k)!}f_k(x_1,\ldots,x_k).
}
While $f_k$ is a probability density, $R_k$ is not.
\end{remark}
We underline that the $0,k$-point and the $k,0$-point correlation functions exist, as we are considering densities on the squared eigenradii and densities on the singular values, and those functions are $k$-point correlation functions~\eqref{eq:k-point function}. To avoid confusion we will refer to $k$-point correlation functions by $0,k$-point or $k,0$-point correlation functions or state clearly whether it refers to $k$ squared singular values or $k$ squared eigenradii.

\nota{}{
\textup{(Probability densities)} We will denote by
\Eq{}{
\rho_{\rm EV} :&\, \mathbb{R}_+ \longrightarrow \mathbb{R}_+, \quad r\mapsto \rho_{\rm EV}(r)=f_{1,0}(r)
} 
the $1$-point correlation function for  the squared eigenradii, and 
\Eq{}{
\rho_{\rm SV} :&\, \mathbb{R}_+ \longrightarrow \mathbb{R}_+,\quad a\mapsto \rho_{\rm SV}(a)=f_{0,1}(a)
} 
the $1$-point correlation function for the squared singular values. 
}

When studying the interaction between singular values and eigenradii, and more generally between two sets of random variables, one also wishes to know and measure how much correlated those random variables are. When the eigenvalues are independent of the singular values, the $j,k$-point correlation function is simply the product of the respective $j$-point and $k$-point correlation functions. Therefore, the difference between the $j,k$-point correlation function and this product quantifies their dependence in the general case. We coin this measure \textit{cross-covariance density}.

\defi{def: cross-covariance}{
\textup{(Cross-covariance density function)} The $j,k$-cross-covariance density function  is defined as the function $cov_{j,k}:\mathbb{R}_+^{j+k}\to \mathbb{R}$ given by
\Eq{}{
\mathrm{cov}_{j,k}(x;y):=f_{j,k}(x;y)-f_{j,0}(x)f_{0,k}(y),
}
where the $f_{j,k}$, $f_{j,0}$ and $f_{0,k}$ are respectively the $j,k$-point, $j,0$-point and $0,k$-point correlation functions as defined in~\eqref{eq:j,k-point pdf}. For convenience, we denote $\mathrm{cov}_{1,1}:=\mathrm{cov}$ and refer to it simply by cross-covariance density function.
}

\begin{remark}
\label{rem:cross-covariance}
 The definition of the cross-covariance is natural. Indeed, when taking two continuous bounded functions $\psi\in \mathrm{C}_b\left(\mathbb{R}_+^j\right)$ and $\varphi \in \mathrm{C}_b\left(\mathbb{R}_+^k\right)$, we have
 \begin{equation}
 \begin{split}
 &\int_{\mathbb{R}_+^{j+k}}  (f_{j,k}(r;a)-f_{j,0}(r)f_{0,k}(a))\psi(r)\varphi(a) dr da\\
 =&\frac{(n-j)!(n-k)!}{(n!)^2}\sum_{\substack{1\leq l_1<\ldots<l_j\leq n\\1\leq m_1<\ldots<m_k\leq n}}\mathrm{Cov}(\psi(r_{l_1},\ldots,r_{l_j}),\varphi(a_{m_1},\ldots,a_{m_k})),
 \end{split}
 \end{equation}
with the covariance
\Eq{}{\mathrm{Cov}(X, Y):=\mathbb{E}\Big[\Big(X-\mathbb{E}[X]\Big) \Big(Y-\mathbb{E}[Y]\Big)\Big].
}

An an example let us consider the case $j=k=1$. When all first and second moments of $r_j$ and $a_k$ exist, it is
\Eq{}{
\int_{\mathbb{R}_+^2}  (f_{1,1}(r;a)-f_{1,0}(r)f_{0,1}(a))ra \,drda=\frac{1}{n^2}\sum_{j,k=1}^n \mathrm{Cov}(r_j, a_k),
}
meaning the average of the covariances between any squared eigenradius with any squared singular value.
Hence, $f_{1,1}-f_{1,0}f_{0,1}$ is the average cross-covariance density between one squared singular value and one squared eigenradius.

We note that this density is very similar to the one when two variables are of the same kind. Then, the covariance density function is defined as
\Eq{}{
\mathrm{cov}(x_1,x_2)=f_2(x_1,x_2)-f_1(x_1)f_1(x_2),
}
involving the $2$-point and $1$-point correlation functions. In this case, the covariance density is the negative of the $2$-level cluster function~\cite{Dyson1962}, which is often used in the physics literature.

\end{remark}

\subsection{Polynomial and P\'olya Ensembles}

In order to find an explicit formula for the $1,k$-point correlation function $f_{1,k}$, one has to use an explicit expression for the density function on the squared singular values $f_{\rm SV}$. A suitable and rather broad class of ensembles are polynomial ensembles~\cite{Kuijlaars2014,Kuijlaars2014a,Kuijlaars2016,Foerster2020} on $ \mathbb{R}_+^n$ which is a probability density function of the form
\Eq{eq:def poly ensemble}{
f_{\rm SV}(x_1,\ldots,x_n)=C_{\rm SV}(w)\Delta_n(x)\det\left[ w_{k-1}(x_j)  \right]_{j,k=1}^n \geq0
}
with
\Eq{}{
\frac{1}{C_{\rm SV}(w)}=n!\,\det\left[\int_{0}^\infty x^{j-1} w_{k-1}(x) dx  \right]_{j,k=1}^n \quad \in \mathbb{R}\setminus \{ 0\}
}
and $\left\{w_{b}\right\}_{b=0}^{n-1}$ , $w_{b} \in \mathrm{L}^1(\mathbb{R}_+)$ whose first $n-1$ moments exist, cf. Eq.~\eqref{eq:polynomial ensemble def}.

A polynomial ensemble is a determinantal point process and its kernel $K$ can be written as follows
\Eq{kernel}{
K(x,y)=\displaystyle\sum_{b=0}^{n-1} q_b(x) p_b(y),
}
with $p_b$ a polynomial of degree $p$ and   $q_{b}\in \mathrm{Span}\{w_0,\dots,w_{n-1}\}$ where $\mathrm{Span}$ is the linear span. The two sets of functions $\{q_b\}_{b=0,\ldots,n-1}$ and  $\{p_b\}_{b=0,\ldots,n-1}$ must satisfy the bi-orthogonality relation
\Eq{biorthogonality}{\int_{0}^\infty q_b(x)p_c(x)dx=\delta_{c,b},}
where $\delta_{c,b}$ is the Kronecker delta function.

 The $1$-point function for a polynomial ensemble, and more generally for any determinantal point process can be expressed in terms of this kernel. It has the rather simple form
\Eq{}{f_1(x)=\frac{1}{n}K(x,x).}

P\'olya ensembles~\cite{Kieburg2015,Kieburg2016,Kieburg2022,Foerster2020} are particular kinds of polynomial ensembles. For those exist a weight function $w$ such that for all $ k \in \llbracket 1,n\rrbracket$, $x\mapsto (-x\partial_x)^{k-1} w(x) \in \mathrm{L}^1(\mathbb{R}_+)$ and that $\mathrm{Span}\{w_{k-1}\}_{k=1}^{n}=\mathrm{Span}\{x\mapsto (-x\partial_x)^{k-1} w(x)\}_{k=1}^{n}$ for which we need an $(n-1)$-times continuous differentiable weight  function $w$. This choice comes very naturally when considering the identity
\Eq{}{
\mathcal{M}[(-x\partial_x)^{n} f(x)](s)=s^n \mathcal{M}f(s), \quad n\in \mathbb{N}
}
for the Mellin transform~\eqref{eq:M-trans} and the fact that the Mellin transform plays an important role in the relation between the joint probability density functions of the singular values and the eigenvalues.

 The joint probability density function of the singular values of a P\'olya ensemble has the form
\Eq{eq:polya ensemble}{
f(x_1,\ldots,x_n)=C_{\rm SV}(w)\Delta_n(x)\det\left[ (-x_j\partial_{x_j})^{k-1} w(x_j)  \right]_{j,k=1}^n .
}
The kernel, then, admits an integral representation, see~\cite[Eq.(4.22)]{Kieburg2016},
\Eq{eq:kernel polya}{
K(x,y)=\sum_{j=0}^{n-1} q_j(x) p_j(y)=n\int_0^1 dt\, q_n(xt)p_{n-1}(yt)
}
with the bi-orthonormal sets of functions $\left\{q_{j}\right\}_{j=0}^{n-1}$ and $\left\{p_{j}\right\}_{j=0}^{n-1}$  given by
\Eq{eq:polya ortho functions}{
p_{j}(x) &=\sum_{c=0}^{j} \binom{j}{c} \frac{(-x)^c}{\mathcal{M}w(c+1)}, \quad q_j(x) = \frac{1}{j!}\partial_x^j[x^j w(x)].}
For this result we actually need the $n$-times differentiability of $w$, i.e., $w\in \mathrm{C}^n(\mathbb{R})$, to have a well-defined function $q_n$. If the weight is only $(n-1)$-differentiable one needs to modify the formula to
\begin{equation}
K(x,y)=n\int_0^1 dt\, q_{n-1}(xt)p_{n-2}(yt)+q_{n-1}(x)p_{n-1}(y).
\end{equation}

Finally, we would like to point out that an orthogonal ensemble~\cite{Anderson2010,Akemann2015,Deift2009,Forrester2010} is a polynomial ensemble for which there exists a weight function $\omega$ such that for all $k \in \llbracket 0,n-1\rrbracket$, $x\mapsto x^k \omega(x) \in \mathrm{L}^1(\mathbb{R}_+)$ and that $\mathrm{Span}\{w_{k-1}\}_{k=1}^{n}=\mathrm{Span}\{x \mapsto x^{k-1}\omega(x)\}_{k=1}^{n}$. The probability density function can therefore be written as follows
\Eq{eq:orthogonal ensemble}{
f(x_1,\ldots,x_n)&=C_{\rm SV}(w)\Delta_n(x)\det\left[ x_j^{k-1} \omega(x_j)  \right]_{j,k=1}^n = C_{\rm SV}(w)\Delta_n(x)^2 \prod_{j=1}^n \omega(x_j) .
}
The orthogonal ensemble is then usually named after the kind of polynomial on which its kernel is built. The classical Laguerre ensemble as well as Jacobi ensemble are both, orthogonal and P\'olya ensembles. Yet, there are P\'olya ensembles that are not orthogonal ensembles and vice versa.

\section{The \textbf{1,k}-point Correlation Function}\label{1,k-point function}

In the ensuing proofs we make use of the following lemma to compute the $n$-fold integral on the eigenangles; cf. \cite[Lemma 1.4]{Kostlan1992}.

\lem{lem:det perm}{
Let $n \in \mathbb{N}$, $n\geq 1$. $r=(r_1,\ldots,r_n)\in \mathbb{R}_+^n$, $\theta=(\theta_1,\ldots,\theta_n)\in [0,2\pi]^{n}$.
\Eq{}{
\int_{[0,2\pi]^{n}} d\theta\ |\Delta_n(\sqrt{r}e^{i\theta})|^2= (2\pi)^n \mathrm{Perm}[r_j^{k-1}]_{j,k=1}^n
}
}

\begin{proof}
    At the heart of the proof lies the integral $\int_0^{2\pi}d\vartheta\ e^{i\vartheta (j-k)}=2\pi \delta_{j,k}$ for all $j,k\in\mathbb{Z}$. Expanding both Vandermonde determinants via the Leibniz formula, we find the claim
\Eq{}{
\int_{[0,2\pi]^{n}} |\Delta_n(\sqrt{r}e^{i\theta})|^2 d\theta &=  \int_{[0,2\pi]^{n}} d\theta \sum_{\substack{\sigma_1,\sigma_2 \in \mathrm{S}_n }} \sign(\sigma_1\sigma_2) \prod_{j=1}^n r_j^{(\sigma_1(j)+\sigma_2(j)-2)/2} e^{i\theta_{j}(\sigma_1(j)-\sigma_2(j))}\\
&= (2\pi)^n \sum_{\sigma_1 \in \mathrm{S}_n} \prod_{j=1}^n r_j^{\sigma_1(j)-1}= (2\pi)^n \mathrm{Perm}[r_j^{k-1}]_{j,k=1}^n ,
}
where $\mathrm{S}_n$ is the symmetric group permuting $n$ elements and ${\rm sign}$ is the signum function which is $+1$ for an even permutation and $-1$ for an odd one.
\end{proof}

\subsection{Proof of \autoref{theo:1,kpt}}\label{Proof of theo:1,kpt}

We will assume $n>2$ in the present section as some steps require this condition to be eligible. Moreover, we assume $f_{\rm SV}\in  \mathrm{L}^{1,\rm SV}(A)$ is a probability density function.

The most important statement of \autoref{theo:1,kpt} is that $f_{1,k}$ is indeed a function and not only a measure. For this purpose we choose a Schwartz function $\phi_{\rm EV}\in \mathcal{S}(\mathbb{R}_+)$ and a continuous bounded function $\phi_{\rm SV}\in \mathrm{C}_{b}(\mathbb{R}_+^{k})$ and make use of~\eqref{relation.eig.sing} in combination with the choice~\eqref{ident.test} for $j=1$, to get
\Eq{}{
&\mathbb{E}\bigg[\frac{1}{n!\,n} \sum_{\substack{\sigma\in \mathrm{S}_n \\1\leq l\leq n }} \phi_{\rm EV}(r_{l}(X)) \phi_{\rm SV}(a_{\sigma(1)}(X),\ldots,a_{\sigma(k)}(X) ) \bigg]\\
=&\, \frac{\left(\prod_{j=0}^{n-1} j!\right)}{\pi^n n(n!)^3} \int_Z dz |\Delta_n(z)|^2\lim_{\varepsilon \to 0} \int_{\mathcal{C}(n)}\left[\prod_{k=1}^n \frac{ds_k}{2\pi i}\zeta(\varepsilon \Im{s_k})\right]\mathrm{Perm}[|z_b|^{-2 s_c}  ]_{b,c=1}^n \\
 & \times \int_A da\, f_{\rm SV}(a)\sum_{\substack{\sigma\in \mathrm{S}_n \\1\leq l\leq n }}\phi_{\rm EV}(|z_l|^2)\phi_{\rm SV}(a_{\sigma(1)},\ldots,a_{\sigma(k)})\frac{\det [a_b^{s_c-1}  ]_{b,c=1}^n }{\Delta_n(a)\Delta_n(s)}.
}
We then integrate over the eigenangles using \autoref{lem:det perm}. Additionally, we can use the permutation symmetry in $r_j$ to replace the sum over $\phi_{\rm EV}(r_l)\phi_{\rm SV}(a_{\sigma(1)},\ldots,a_{\sigma(k)})$ by $n\phi_{\rm EV}(r_1)\phi_{\rm SV}(a_{\sigma(1)},\ldots,a_{\sigma(k)})$. We are left with
\Eq{}{
&\mathbb{E}\bigg[\frac{1}{n!\,n} \sum_{\substack{\sigma\in \mathrm{S}_n \\1\leq l\leq n }} \phi_{\rm EV}(r_{l}(X)) \phi_{\rm SV}(a_{\sigma(1)}(X),\ldots,a_{\sigma(k)}(X) ) \bigg]\\
=&\, \frac{\left(\prod_{j=0}^{n-1} j!\right)}{ (n!)^3} \int_A dr\,\mathrm{Perm}[r_b^{c-1}  ]_{b,c=1}^n\lim_{\varepsilon \to 0} \int_{\mathcal{C}(n)}\left[\prod_{k=1}^n \frac{ds_k}{2\pi i}\zeta(\varepsilon \Im{s_k})\right] \\
 & \times \mathrm{Perm}[r_b^{- s_c}  ]_{b,c=1}^n\int_A da\, f_{\rm SV}(a)\sum_{\sigma\in \mathrm{S}_n}\phi_{\rm EV}(r_1)\phi_{\rm SV}(a_{\sigma(1)},\ldots,a_{\sigma(k)})\frac{\det [a_b^{s_c-1}  ]_{b,c=1}^n }{\Delta_n(a)\Delta_n(s)}.
}
We expand the first permanent in the first column and integrate first over $r_2,\ldots,r_n$ and then over $r_1$. Particularly the permutation invariance of the remaining integrand in $r_2,\ldots,r_n$ tells us that each integral has $(n-1)!$ identical contributions so that
\Eq{}{
&\mathbb{E}\bigg[\frac{1}{n!\,n} \sum_{\substack{\sigma\in \mathrm{S}_n \\1\leq l\leq n }} \phi_{\rm EV}(r_{l}(X)) \phi_{\rm SV}(a_{\sigma(1)}(X),\ldots,a_{\sigma(k)}(X) ) \bigg]\\
=&\, \frac{\left(\prod_{j=0}^{n-1} j!\right)}{ (n!)^2 n} \sum_{l=1}^{n}\int_0^\infty dr_1\, r_1^{l-1}\left[\prod_{j=2}^ndr_j r_j^{\tau_{j}(l)-1}\right] \lim_{\varepsilon \to 0}\int_{\mathcal{C}(n)}\left[\prod_{k=1}^n \frac{ds_k}{2\pi i}\zeta(\varepsilon \Im{s_k})\right]\\
 & \times \mathrm{Perm}[r_b^{- s_c}  ]_{b,c=1}^n  \int_A da\, f_{\rm SV}(a)\sum_{\sigma\in \mathrm{S}_n}\phi_{\rm EV}(r_1)\phi_{\rm SV}(a_{\sigma(1)},\ldots,a_{\sigma(k)})\frac{\det [a_b^{s_c-1}  ]_{b,c=1}^n }{\Delta_n(a)\Delta_n(s)}
}
with the $(n-1)$-dimensional vector $\tau(l)=(1,\ldots,l-1,l+1,\ldots,n)$.

We define the function
\begin{equation}
g(a):=\frac{1}{n!\,n} f_{\rm SV}(a)\sum_{\sigma\in \mathrm{S}_n}\phi_{\rm SV}(a_{\sigma(1)},\ldots,a_{\sigma(k)})
\end{equation}
which is in $\mathrm{L}^{1,\rm SV}(A)$. Now, we can identify the integral over $a$ with the spherical transform~\eqref{eq:S-trans} of $g$, the integral over $s$ and the limit $\varepsilon\to0$ with the multivariate inverse Mellin transform~\eqref{inv.Mellin} and the integral over $r_2,\ldots,r_n$ as the Mellin transform ${\rm id}\otimes\mathcal{M}^{n-1}$ in the last $n-1$ entries, i.e.,
\Eq{eq: 3.7}{
&\mathbb{E}\bigg[\frac{1}{n!\,n} \sum_{\substack{\sigma\in \mathrm{S}_n \\1\leq l\leq n }} \phi_{\rm EV}(r_{l}(X)) \phi_{\rm SV}(a_{\sigma(1)}(X),\ldots,a_{\sigma(k)}(X) ) \bigg]\\
=&\, \left(\prod_{j=0}^{n-1} j!\right) \sum_{l=1}^{n}\int_0^\infty dr_1\, r_1^{l-1}\phi_{\rm EV}(r_1) ( {\rm id}\otimes\mathcal{M}^{n-1})\mathcal{M}_S^{-1}\mathcal{S}g(r_1,\tau(l))\\
=& \left(\prod_{j=0}^{n-1} j!\right) \sum_{l=1}^{n}\int_0^\infty dr_1\, r_1^{l-1} \phi_{\rm EV}(r_1)( \mathcal{M}^{-1}\otimes{\rm id}^{\otimes n-1})\mathcal{S}g(r_1,\tau(l)).
}
Note that $r_1\mapsto r_1^{l-1} (\mathcal{M}^{-1}\otimes{\rm id}^{\otimes n-1})\mathcal{S}g(r_1,\tau(l))\in \mathrm{L}^1(\mathbb{R}_+)$ due to \autoref{theo:R-transform}, as $g\in\mathrm{L}^{1,\rm SV}(A)$. When writing the RHS of \eqref{eq: 3.7} explicitly we arrive at
\Eq{eq:3.8}{
&\mathbb{E}\bigg[\frac{1}{n!\,n} \sum_{\substack{\sigma\in \mathrm{S}_n \\1\leq l\leq n }} \phi_{\rm EV}(r_{l}(X)) \phi_{\rm SV}(a_{\sigma(1)}(X),\ldots,a_{\sigma(k)}(X) ) \bigg]\\
=&\, \frac{\left(\prod_{j=0}^{n-1} j!\right)}{ n!\, n}  \sum_{l=1}^{n}\int_0^\infty dr_1\, r_1^{l-1} \phi_{\rm EV}(r_1)\lim_{\varepsilon\to 0}\int_{\mathcal{C}_l}\frac{ds_l}{2i\pi}\zeta(\varepsilon \Im{s_l}) r_1^{-s_l}\\
&\times\int_A da\, f_{\rm SV}(a)\sum_{\sigma\in \mathrm{S}_n}\phi_{\rm SV}(a_{\sigma(1)},\ldots,a_{\sigma(k)})\frac{\det [a_b^{s_l-1},\ a_b^{\tau_c(l)-1}  ]_{\substack{b=1,\ldots,n\\ c=1,\ldots,n-1}} }{\Delta_n(a)\Delta_{n-1}(\tau(l))\prod_{j=1}^{n-1}(\tau_j(l)-s_l)},
}
where we recall that $\mathcal{C}_l= l+i\mathbb{R}$ and where we have used the fact
\Eq{}{
\Delta_{n}(s_l,\tau(l))=\Delta_{n-1}(\tau(l))\prod_{j=1}^{n-1}(\tau_j(l)-s_l).
}
To carry out the limit $\varepsilon\to0$ ones needs to be aware that the integral over $s_l$ is not necessarily absolutely convergent, namely in the case of degenerate $a$. We will see that the following order of integration
\Eq{}{
\lim_{\varepsilon\to 0}\int da \int d r_1 \int d s_l.
}
enables us to carry out the limit explicitly with the help of Lebesgue's dominated convergence theorem. To achieve this order, we first need to show that the limit $\varepsilon\to0$ can be interchanged with the integration over $r_1$. The Schwartz function $\phi_{\rm EV}$ and the regularisation function $\zeta$ will help us achieve it. The idea is to rewrite the complex integral in the RHS of \eqref{eq:3.8} in the form of a multiplicative convolution. Indeed, the complex integral over $s_l$ itself (without the limit $\varepsilon\to 0$) corresponds to the usual (unregularized) inverse Mellin transform of the product of the two functions
\Eq{}{
\tilde{\zeta}_\eps(s):=\zeta(\varepsilon \Im{s}),\quad \tilde{g}_l(s):= \mathcal{S}g(s,\tau(l)).
}
For the sake of clarity, we will denote by $\mathcal{M}^{-1}_*$ the unregularised inverse Mellin transform. The integrals in the RHS of \eqref{eq:3.8} can then be written in the compact form
\Eq{}{
\mathcal{M}\left[\phi_{\rm EV}\cdot\lim_{\eps\to0}\mathcal{M}^{-1}_*[\tilde{\zeta}_\eps\cdot \tilde{g}_l]\right](l).
}
Using Mellin inversion theorem, the function $\mathcal{M}^{-1}_*[\tilde{\zeta}_\eps\cdot \tilde{g}_l]$ can expressed in terms of a Mellin convolution
\Eq{eq: inv mellin convo}{
\int_{\mathcal{C}_l}\frac{ds_l}{2i\pi} r_1^{-s_l} \tilde{\zeta}_\eps(s_l) \tilde{g}_l(s_l)=\mathcal{M}^{-1}[\tilde{\zeta}_\eps\cdot \tilde{g}_l](r_1)=\int_0^\infty \frac{dy}{y} \mathcal{M}^{-1}_*\tilde{\zeta}_\eps(y) \mathcal{M}^{-1}\tilde{g}_l\left(\frac{r_1}{y} \right),
}
where the inverse Mellin transform of $\tilde{g}_l$ remains regularised  and implicit, i.e. $\mathcal{M}^{-1}\tilde{g}_l(y)=(\mathcal{M}^{-1}\otimes{\rm id}^{\otimes n-1})\mathcal{S}g(y,\tau(l))$, and the inverse Mellin transform of the auxiliary function $\tilde{\zeta}_\eps$, see~\eqref{eq:reg fct}), is 
\Eq{eq: inv mellin zeta}{
\mathcal{M}^{-1}_*\tilde{\zeta}_\eps\left(x \right) =\int_{\mathcal{C}_l}\frac{ds_l}{2i\pi}\zeta(\varepsilon \Im{s_l}) x^{-s_l}=\frac{\pi}{4\varepsilon}x^{-l}\cos\left(\frac{\pi \log(x)}{2\varepsilon}\right)\left(\Theta\left(e^\varepsilon-x\right)-\Theta\left(e^{-\varepsilon}-x\right) \right),
}
with $\Theta$ the Heaviside step function. We then have implicitly added another regularisation through $\mathcal{M}^{-1}\tilde{g}_l$.

Thus, combining \eqref{eq: inv mellin convo} and \eqref{eq: inv mellin zeta}, and applying already the Heaviside step functions, we have
\Eq{}{
&\mathcal{M}\left[\phi_{\rm EV}\cdot\mathcal{M}^{-1}_*[\tilde{\zeta}_\eps\cdot \tilde{g}_l]\right](l)\\
=&\int_0^\infty dr_1\, r_1^{l-1} \phi_{\rm EV}(r_1)\int_{e^{-\eps}}^{e^\eps}\frac{dy}{y}\frac{\pi}{4\varepsilon}y^{-l}\cos\left(\frac{\pi \log(y)}{2\varepsilon}\right)\mathcal{M}^{-1}\tilde{g}_l\left(\frac{r_1}{y} \right).
}
Then, one can interchange the two integrals. This is allowed  due to the integrability of  $\phi_{\rm EV}$ over the compact set of a given fixed $\varepsilon>0$ and that $r_1\mapsto r_1^{l-1}\mathcal{M}^{-1}\tilde{g}_l(r_1)\in\mathrm{L}^1(\mathbb{R}_+)$, as mentioned earlier. Performing the consecutive changes of variable $r_1\mapsto r_1 y$ and $y=e^{\eps x}$ yields
\Eq{eq: 3.15}{
&\mathcal{M}\left[\phi_{\rm EV}\cdot\mathcal{M}^{-1}_*[\tilde{\zeta}_\eps\cdot \tilde{g}_l]\right](l)\\
=&\int_0^\infty dr_1\, r_1^{l-1}\left[ \int_{-1}^{1}dx\, \phi_{\rm EV}(e^{\varepsilon x}r_1)\frac{\pi}{4}\cos\left(\frac{\pi x}{2}\right)\right](\mathcal{M}^{-1}\otimes{\rm id}^{\otimes n-1})\mathcal{S}g(r_1,\tau(l)),
}
which can be compared to the integrals \eqref{eq:3.8}, where the only difference is the absence of $\lim_{\eps\to 0}$.  Taking now the limit $\eps\to 0$ on the outside of \eqref{eq: 3.15}, it is now clear that one can use Lebesgue's dominated convergence theorem to exchange $\lim_{\eps\to 0}$ with the integral over $r_1$. Indeed, the integral over $y$ is bounded by a constant independent of $r_1$ and $\eps$ and has a pointwise limit as $\eps \to 0$ because $\phi_{\rm EV}$ is a Schwartz function. The integrability in $r_1$ is guaranteed by $r_1\mapsto r_1^{l-1}\mathcal{M}^{-1}\tilde{g}_l(r_1)\in\mathrm{L}^1(\mathbb{R}_+)$. Finally,
\Eq{}{
\lim_{\eps\to 0}\mathcal{M}\left[\phi_{\rm EV}\cdot\mathcal{M}^{-1}_*[\tilde{\zeta}_\eps\cdot \tilde{g}_l]\right](l)=\mathcal{M}\left[\phi_{\rm EV}\cdot \lim_{\eps\to 0}\mathcal{M}^{-1}_*[\tilde{\zeta}_\eps\cdot \tilde{g}_l]\right](l).
}

Thus, keeping $\lim_{\eps\to 0}$ outside of all the integrals (LHS of the above equation) and going on with our original aim of proving \autoref{theo:1,kpt}, we can now interchange the integrals over $r_1$, $s_l$ and $a$ as we like since the absolute integrability in $s_l$ is given by the auxiliary function $\zeta$, the one in $r_1$ by $\phi_{\rm EV}$ and the one in $a$ by $f_{\rm SV}$. Setting $s_l=l+it$, with $t\in\mathbb{R}$, we note that the last term in \eqref{eq:3.8} is not a problem because 
\begin{equation}\label{eq: GN int}
\begin{split}
&\left|\frac{\left(\prod_{j=0}^{n-1} j!\right)\det [a_b^{l+it-1},\ a_b^{\tau_c(l)-1}  ]_{\substack{b=1,\ldots,n\\ c=1,\ldots,n-1}} }{\Delta_n(a)\Delta_{n-1}(\tau(l))\prod_{j=1}^{n-1}(\tau_j(l)-l-it)}\right|\\
=&\left|\int_{{\rm U}(n)}\left(\frac{\det\left( (U\diag(a_1,\ldots,a_n)U^\dagger)_{(n-l)\times (n-l)}\right)}{\det \left((U\diag(a_1,\ldots,a_n)U^\dagger)_{(n-l-1)\times (n-l-1)}\right)}\right)^{it}d\mu_H(U)\right|\leq1,
\end{split}
\end{equation}
where $(U\diag(a_1,\ldots,a_n)U^\dagger)_{p\times p}$ is the upper left $p\times p$ block and $d\mu_H(U)$ is the normalised Haar measure of the unitary group. The group integral in the second line  is a particular case of the Gelfand-Naimark integral~\cite{Gelfand1950} which evidently is an integral of a bounded integrand over the unitary group; cf. \cite[Subsection 2.4]{Kieburg2016}.

Having interchanged the limit $\eps\to 0$ with the $r_1$ integral, we then interchange the $t$ and $r_1$ integrals with the $a$ integral to slightly shift the contour of $t\to t+i\delta$ with a $\delta>0$, i.e.,
\Eq{}{
&\mathbb{E}\bigg[\frac{1}{n!\,n} \sum_{\substack{\sigma\in \mathrm{S}_n \\1\leq l\leq n }} \phi_{\rm EV}(r_{l}(X)) \phi_{\rm SV}(a_{\sigma(1)}(X),\ldots,a_{\sigma(k)}(X) ) \bigg]\\
=&\, \frac{\left(\prod_{j=0}^{n-1} j!\right)}{ n!\, n}  \sum_{l=1}^{n}\lim_{\varepsilon\to0}\int_A da\, f_{\rm SV}(a)\sum_{\sigma\in \mathrm{S}_n}\phi_{\rm SV}(a_{\sigma(1)},\ldots,a_{\sigma(k)})\int_0^\infty dr_1\, \phi_{\rm EV}(r_1) \\
&\times\int_{-\infty}^{\infty}\frac{dt}{2\pi}\zeta(\varepsilon (t+i\delta)) r_1^{-1+\delta-it}\frac{\det [a_b^{l-1-\delta+it},\ a_b^{\tau_c(l)-1}  ]_{\substack{b=1,\ldots,n\\ c=1,\ldots,n-1}} }{\Delta_n(a)\Delta_{n-1}(\tau(l))\prod_{j=1}^{n-1}(\tau_j(l)-l+\delta-it)}.
}
The shift is allowed as the integrand in $t$ is an entire function and it drops off in the real direction at least like $t^{-2}$ because of $\zeta$.

The shift allows us to interchange the $r_1$ and $t$ integrals so that we can integrate by parts twice in  $r_1$,
\begin{equation}
\int_0^\infty dr_1\, \phi_{\rm EV}(r_1)r_1^{-1+\delta-it}=\frac{1}{(\delta -it)(1+\delta-it)}\int_0^\infty dr_1\, \phi_{\rm EV}''(r_1)r_1^{1+\delta-it}
\end{equation}
which creates an additional decrease in $t$ for $t\to0$. That decrease will now take over the job of the auxiliary function $\zeta$. Since the integrability in $a$ might be an issue because of the shift, we take the limit $\delta\to0$ which gives $-1/2$ times the residue at $t=-i\delta=0$ and a Cauchy principal value integral, denoted $\dashint$, with a simple pole at $t=0$. This leads to
\Eq{eq: 3.19}{
&\mathbb{E}\bigg[\frac{1}{n!\,n} \sum_{\substack{\sigma\in \mathrm{S}_n \\1\leq l\leq n }} \phi_{\rm EV}(r_{l}(X)) \phi_{\rm SV}(a_{\sigma(1)}(X),\ldots,a_{\sigma(k)}(X) ) \bigg]\\
=&\, \frac{\left(\prod_{j=0}^{n-1} j!\right)}{ n!\, n}  \sum_{l=1}^{n}\lim_{\varepsilon\to0}\int_A da\, f_{\rm SV}(a)\sum_{\sigma\in \mathrm{S}_n}\phi_{\rm SV}(a_{\sigma(1)},\ldots,a_{\sigma(k)})\int_0^\infty dr_1\, \phi_{\rm EV}''(r_1) r_1\\
&\times\left[\frac{1}{2(\prod_{j=0}^{n-1}j!)}-\mathcal{I}\right],
}
with the Cauchy principal value integral
\Eq{}{
\mathcal{I}:=\dashint_{-\infty}^{\infty}\frac{dt}{2\pi i t(1-it)}\zeta(\varepsilon t) r_1^{-it}\frac{\det [a_b^{l+it-1},\ a_b^{\tau_c(l)-1}  ]_{\substack{b=1,\ldots,n\\ c=1,\ldots,n-1}} }{\Delta_n(a)\Delta_{n-1}(\tau(l))\prod_{j=1}^{n-1}(\tau_j(l)-l-it)}.
}
We can get back to an ordinary integral without principal value when defining the function
\begin{equation}
\Upsilon_l(t;r_1;a):=\frac{1}{ t\Delta_n(a)\Delta_{n-1}(\tau(l))}\Im{\frac{\det [ r_1^{-it}a_b^{l+it-1},\ a_b^{\tau_c(l)-1}  ]_{\substack{b=1,\ldots,n\\ c=1,\ldots,n-1}} }{(1-it)\prod_{j=1}^{n-1}(\tau_j(l)-l-it)}}
\end{equation}
and note that, as $\mathcal{I}$ is real, and due to the symmetry $\zeta(\varepsilon t)=\zeta(-\varepsilon t)$, one has
\begin{equation}
\mathcal{I}=\frac{1}{2}(\mathcal{I}+\mathcal{I}^*)=\int_{-\infty}^{\infty}\frac{dt}{2\pi }\zeta(\varepsilon t) \Upsilon_l(t;r_1;a),
\end{equation}
where we have used the fact $\Re{-iz}=\Im{z}$ and denoted by $\mathcal{I}^*$ the complex conjugate of $\mathcal{I}$. The function $\Upsilon_l(t;r_1;a)$ has no pole at $t=0$ and is absolutely integrable in $t$ as it behaves like $O(1/t^2)$ for $t\to\pm\infty$ for all $a\in\mathbb{R}_+^k$ and $r_1\in\mathbb{R}_+$.

The full integrand, in the RHS of \eqref{eq: 3.19}, is now bounded from above by the integrable function $f_{\rm SV}(a)|r_1\phi_{\rm EV}''(r_1)\Upsilon_l(t;r_1;a)|$ independent of $\varepsilon$ since the rest is bounded by a constant while the pointwise limit $\varepsilon\to0$ is clear. Thence, we can carry out this limit via Lebesgue's dominated convergence theorem and arrive at
\Eq{}{
&\mathbb{E}\bigg[\frac{1}{n!\,n} \sum_{\substack{\sigma\in \mathrm{S}_n \\1\leq l\leq n }} \phi_{\rm EV}(r_{l}(X)) \phi_{\rm SV}(a_{\sigma(1)}(X),\ldots,a_{\sigma(k)}(X) ) \bigg]\\
=&\int_A da\, f_{\rm SV}(a)\phi_{\rm SV}(a_{1},\ldots,a_{k})\int_0^\infty dr_1\, \phi_{\rm EV}''(r_1) r_1\biggl[\frac{1}{2}-\frac{\prod_{j=0}^{n-1} j!}{n}\sum_{l=1}^{n}\int_{-\infty}^{\infty}\frac{dt}{2\pi }\Upsilon_l(t;r_1;a)\biggl],
}
where we have made use of the permutation invariance of the integrand in $a$ which yields an $n!$ factor.
To figure out the analytic properties of the function in the bracket it is helpful to plug in the Gelfand-Naimark integral~\eqref{eq: GN int} and interchange the integral over the unitary matrix $U$ with the $t$-integral which can be carried out with standard residue theorem techniques yielding
\Eq{}{
&r_1\biggl[\frac{1}{2}-\frac{\left(\prod_{j=0}^{n-1} j!\right)}{n}\sum_{l=1}^{n}\int_{-\infty}^{\infty}\frac{dt}{2\pi }\Upsilon_l(t;r_1;a)\biggl]\\
=&\frac{1}{n}\sum_{l=1}^{n}\int_{{\rm U}(n)}d\mu_H(U)\left[r_1-\frac{\det\left( (U\diag(a_1,\ldots,a_n)U^\dagger)_{(n-l)\times (n-l)}\right)}{\det \left((U\diag(a_1,\ldots,a_n)U^\dagger)_{(n-l-1)\times (n-l-1)}\right)}\right]\\
&\times\Theta\left(r_1-\frac{\det\left( (U\diag(a_1,\ldots,a_n)U^\dagger)_{(n-l)\times (n-l)}\right)}{\det \left((U\diag(a_1,\ldots,a_n)U^\dagger)_{(n-l-1)\times (n-l-1)}\right)}\right).
}
This function is evidently once piecewise differentiable and continuous so we can proceed with an integration by parts,
\Eq{}{
&\mathbb{E}\bigg[\frac{1}{n!\,n} \sum_{\substack{\sigma\in \mathrm{S}_n \\1\leq l\leq n }} \phi_{\rm EV}(r_{l}(X)) \phi_{\rm SV}(a_{\sigma(1)}(X),\ldots,a_{\sigma(k)}(X) ) \bigg]\\
=&-\int_A da\, \int_0^\infty dr_1\, f_{\rm SV}(a)\phi_{\rm SV}(a_{1},\ldots,a_{k})\phi_{\rm EV}'(r_1) \Xi(r_1;a)
}
with
\begin{equation}\label{Xi.def}
\begin{split}
\Xi(r_1;a):=&\frac{1}{n}\sum_{l=1}^{n}\int_{{\rm U}(n)}d\mu_H(U)\Theta\left(r_1-\frac{\det\left( (U\diag(a_1,\ldots,a_n)U^\dagger)_{(n-l)\times (n-l)}\right)}{\det \left((U\diag(a_1,\ldots,a_n)U^\dagger)_{(n-l-1)\times (n-l-1)}\right)}\right)\\
=&\frac{\left(\prod_{j=0}^{n-1} j!\right)}{n}\sum_{l=1}^{n}\int_{-\infty}^{\infty}\frac{dt\,r_1^{-it+1/2}}{2\pi (i t-1/2)}\frac{\det [ a_b^{l-3/2+it},\ a_b^{\tau_c(l)-1}  ]_{\substack{b=1,\ldots,n\\ c=1,\ldots,n-1}} }{\Delta_n(a)\Delta_{n-1}(\tau(l))\prod_{j=1}^{n-1}(\tau_j(l)-l+1/2-it)}.
\end{split}
\end{equation}
The first expression, in terms of a group integral, makes it immediate that it does not cause any integrability issues in $a$ and $r_1$.
The latter integral has to be understood in the sense of Cauchy's principal value integral if all squared singular values $a$ are equal to each other, and is absolutely convergent otherwise.

We are now ready to finish the proof. The function $ f_{\rm SV}$ is a density with respect to the Lebesgue measure and the event of two singular values being equal to each other has Lebesgue measure $0$. Hence, we can consider the $r_1$ integral in terms of fixed pairwise distinct square singular values. This tells us that the integrand in~\eqref{Xi.def} behaves like $1/t^n$ for $t\to\infty$ and vanishes at $r_1=0$. For $n>2$ this is continuously differentiable which allows us to integrate by parts once more, i.e.,
\Eq{}{
&\mathbb{E}\bigg[\frac{1}{n!\,n} \sum_{\substack{\sigma\in \mathrm{S}_n \\1\leq l\leq n }} \phi_{\rm EV}(r_{l}(X)) \phi_{\rm SV}(a_{\sigma(1)}(X),\ldots,a_{\sigma(k)}(X) ) \bigg]\\
=&\int_A da\, \int_0^\infty dr_1\, f_{\rm SV}(a)\phi_{\rm SV}(a_{1},\ldots,a_{k})\phi_{\rm EV}(r_1) \Xi'(r_1;a)
}
with
\begin{equation}\label{Xi.prime}
\begin{split}
\Xi'(r_1;a)=&\frac{\left(\prod_{j=0}^{n-1} j!\right)}{n}\sum_{l=1}^{n}\int_{-\infty}^{\infty}\frac{dt}{2\pi }r_1^{-it-1/2}\frac{\det [ a_b^{l-3/2+it},\ a_b^{\tau_c(l)-1}  ]_{\substack{b=1,\ldots,n\\ c=1,\ldots,n-1}} }{\Delta_n(a)\Delta_{n-1}(\tau(l))\prod_{j=1}^{n-1}(\tau_j(l)-l+1/2-it)}.
\end{split}
\end{equation}
The last step is to set back $s=l+it$ and shift $s\mapsto s+1/2$ such that $s\in\mathcal{C}_l=l+i\mathbb{R}$, which is allowed due to the integrand being entire and dropping off rapidly enough. For $k=n>2$, this immediately yields the claim of the theorem. 
For $k<n$ we only integrate out the variables $a_{k+1},\ldots, a_n$. This finishes the proof of the main statement of \autoref{theo:1,kpt}.

We note that the statement of the limit to  degenerate $a_1,\ldots,a_k$ is immediate as the integral over $t$ converges at least conditionally as long as not all singular values are equal to each other. In the latter case when $a_1=\ldots=a_n$, the statement of the weak convergence becomes also clear as we have
\begin{equation}
\Xi(r_1;a)=\Theta(r_1-a_1)
\end{equation}
which implies
\begin{equation}
-\int_0^\infty dr_1\, \phi_{\rm EV}'(r_1) \Xi(r_1;a)=\phi_{\rm EV}(a_1),
\end{equation}
meaning the $1,n$-point measure comprises a Dirac delta function $\delta(r_1-a_1)$ and does not have a density anymore.

\subsection{Proof of \autoref{cor:cond.dens}} \label{sec:cor.proof.1.1}

We start with the result~\eqref{eq:1,1pt} for $k=n$, i.e.,
\begin{equation}
\begin{split}
f_{1,k}(r;a_1,\ldots,a_n)=& \, \frac{1}{n} \left(\prod_{p=0}^{n-1} p! \right) \sum_{j=1}^{n} \int_{\mathcal{C}_j}\frac{ds}{2\pi i}  r^{j-1-s}   f_{\rm SV}(a) \frac{
\det \begin{bmatrix}
    a_b^{s-1} \\
    a_b^{\tau_c(j)-1}\\
  \end{bmatrix}_{\substack{b=1,\ldots,n\\ c=1,\ldots,n-1}} }{\Delta_n(s,\tau(j))\Delta_n(a)}.
\end{split}
\end{equation}
As we do not carry out any integration over $a$ the factor $f_{\rm SV}(a)$ is a constant in the remaining integration variable $s$ and summing index $j$. Thus, the division by $f_{\rm SV}(a)$ immediately yields the first line of~\eqref{cond.dens}.

To get the second result, we  first note that the integrand is entire in $s$, as the apparent $(n-1)$ simple poles cancel with the zeros generated by the determinant in the numerator. Therefore, we can shift the contour $\mathcal{C}_j$ to  the imaginary axis $\mathcal{C}_{0}=i\mathbb{R}$. Next, we evaluate the Vandermode determinant
\Eq{Vand.tau}{\frac{1}{\Delta_{n}(s,\tau(j))}=\frac{(j-1)!\,(n-j)!}{\prod_{l=1}^{n-1} l!} \frac{j-s}{\prod_{l=1}^n(l-s)}.
}
and rewrite $(j-s)r^{j-1-s}=\partial_r r^{j-s}$. The derivative can be swapped with the integral over $s$ because of the absolute integrability in $s$ and holomorphicity in $r$ and $s$ (along $\mathcal{C}_{0}$)  of  the integrand. This means we consider
\begin{equation}
\rho_{\rm EV}(r|a)=\frac{1}{n}\partial_r\sum_{j=1}^{n} (j-1)!\,(n-j)!\int_{\mathcal{C}_{0}}\frac{ds}{2\pi i}     \frac{
\det \begin{bmatrix}
    (a_b/r)^{s-1} \\
    (a_b/r)^{\tau_c(j)-1}\\
  \end{bmatrix}_{\substack{b=1,\ldots,n\\ c=1,\ldots,n-1}} }{\Delta_n(a/r) \prod_{l=1}^n(l-s) },
\end{equation}
where $a/r$ is a shorthand notation for the vector $(a_1/r,\ldots,a_n/r)$. The sum over $j$ can be understood as a Laplace expansion in the first column of an $(n+1)\times(n+1)$ determinant, while the integral over $s$ can be pulled into the first row which can then be computed by closing the contour either around the positive real line for $r>a_b$ or the negative real line for $a_b>r$, i.e.,
\begin{equation}
\begin{split}
\int_{\mathcal{C}_{0}}\frac{ds}{2\pi i}     \frac{ (a_b/r)^{s-1}}{\prod_{l=1}^n(l-s)}=&\Theta(r-a_b)\sum_{l=1}^{n}\frac{(-a_b/r)^{l-1}}{(l-1)!(n-l)!}=\frac{(1-a_b/r)^{n-1}}{(n-1)!}\Theta(r-a_b)
\end{split}
\end{equation}
with $\Theta$ the Heaviside step function. We  arrive at
\begin{equation}
\begin{split}
\rho_{\rm EV}(r|a)=&\frac{1}{n}\partial_r\frac{
\det \left[\begin{array}{c|c} 0 &
    (1-a_b/r)^{n-1}\Theta(r-a_b)  \\\hline \displaystyle (-1)^{c}\frac{(c-1)!(n-c)!}{(n-1)!} &
    (a_b/r)^{c-1} \\
  \end{array}\right]_{b,c=1}^{n} }{\Delta_n(a/r)}.
  \end{split}
\end{equation}
After pulling the factors $(-1)^{c-1}\frac{(c-1)!(n-c)!}{(n-1)!}$ out of the numerator and combining those with the denominator we find the second line of~\eqref{cond.dens}.

The claims about the degeneracy are again clear when taking the limits, in particular that the fully degenerate case can be only understood in the weak sense.

\subsection{Proof of \autoref{theo:poly ensemble}} \label{sec:proof.1.3}

To prove \autoref{theo:poly ensemble}, we start from the following expression of the conditional level density of the squared radii
\begin{equation}
\begin{split}
\rho_{\rm EV}(r|a)=\frac{
\partial_r \det \left[\begin{array}{c|c} 0 &
   \displaystyle -\left(1-\frac{a_b}{r}\right)^{n-1}\Theta(r-a_b)  \\\hline \displaystyle \frac{(c-1)!(n-c)!}{n!}  (-r)^{c-1} &
   \displaystyle a_b^{c-1} 
  \end{array}\right]_{b,c=1}^{n} }{\Delta_n(a)},
  \end{split}
\end{equation}
where $b$ is a column index and $c$ a row index. This equation is a direct consequence of~\eqref{cond.dens}. 

We integrate against the joint probability density  of the squared singular values $f_{\rm SV}$ in $a_{k+1},\ldots,a_n$ which follows from a polynomial ensemble, i.e., $f_{\rm SV}$ is given by~\eqref{eq:polynomial ensemble def}. The Vandermonde determinant $\Delta_n(a)$ cancels so that we have for the $1,k$-point correlation function
\begin{equation}
\begin{split}
f_{1,k}(r;a_1,\ldots,a_k)=&\int_{\mathbb{R}_+^{n-k}}da_{k+1}\cdots da_n\frac{ \det\left[ w_{k-1}(a_j)  \right]_{j,k=1}^n}{n!\,\det\left[ \mathcal{M}w_{k-1}(j)  \right]_{j,k=1}^n}\\
&\times\partial_r \det \left[\begin{array}{c|c} 0 &
   \displaystyle -\left(1-\frac{a_b}{r}\right)^{n-1}\Theta(r-a_b)  \\\hline  \displaystyle \frac{(c-1)!(n-c)!}{n!}  (-r)^{c-1} &
   \displaystyle a_b^{c-1} 
  \end{array}\right]_{b,c=1}^{n} .
\end{split}
\end{equation}
Since $n>2$ the integrand as well as its derivative in $r$ are continuous and integrable due to the weight functions. Thence, we can pull the derivative in $r$ in front of the integral. Furthermore, we can rewrite
\begin{equation}
\frac{(c-1)!(n-c)!}{n!}  (-r)^{c-1}=\int_0^\infty \frac{dt}{(1+t)^{n+1}}(-rt)^{c-1},
\end{equation}
which allows us to recombine the monomials to the monic polynomials $\{p_j\}_{j=0,\ldots,n-1}$ which are biorthogonal to the weights $\{q_j\}_{j=0,\ldots,n-1}$ that span the vector space ${\rm span}\{w_0,\ldots,w_n\}$, see~\eqref{biorthogonality}. This means that we can compute the integral over $a_{k+1},\ldots,a_n$,
\begin{equation}
\begin{split}
f_{1,k}(r;a_1,\ldots,a_k)=&\frac{1}{n!}\partial_r \int_{\mathbb{R}_+^{n-k}}da_{k+1}\cdots da_n\int_0^\infty \frac{dt}{(1+t)^{n+1}} \det\left[ q_{d-1}(a_e)  \right]_{d,e=1}^n\\
&\times\det \left[\begin{array}{c|c} 0 &
   \displaystyle -\left(1-\frac{a_b}{r}\right)^{n-1}\Theta(r-a_b)  \\\hline  \displaystyle p_{c-1}(-rt) &
   \displaystyle p_{c-1}(a_b) 
  \end{array}\right]_{b,c=1}^{n} ,
\end{split}
\end{equation}
with the help of a generalisation~\cite[Appendix C.1]{Kieburg2010} of the Andr\'eief identity~\cite{Andreev1886}. We are allowed to interchange the two sets of integrals due to the absolute integrability. We arrive at
\begin{equation}
\begin{split}
f_{1,k}(r;a_1,\ldots,a_k)=&\frac{(n-k)!}{n!}\partial_r\int_0^\infty \frac{dt}{(1+t)^{n+1}}\\
&\hspace*{-2cm}\times\det \left[\begin{array}{c|c|c} 0 & 0 & q_{d-1}(a_e) \\\hline
 0 &
   \displaystyle -\left(1-\frac{a_b}{r}\right)^{n-1}\Theta(r-a_b) &  \displaystyle-\overset{\ }{\int_0^r} dv\ q_{d-1}(v)\left(1-\frac{v}{r}\right)^{n-1}  \\\hline  \displaystyle p_{c-1}(-rt) &
   \displaystyle p_{c-1}(a_b) & \delta_{cd}
  \end{array}\right],
\end{split}
\end{equation}
where $b=1,\ldots,k$ and $d=1,\ldots,n$ are column indices and $c=1,\ldots,n$ and $e=1,\ldots,k$ are row indices.
We have employed the bi-orthogonality~\eqref{biorthogonality}.

Next, we exploit the standard identity
\Eq{}{\det\left(\begin{array}{c c} 
    	A & B\\ 
    	C & D\\
\end{array}\right)=\det(D)\det \left(A-BD^{-1}C\right)}
and identify the kernel~\eqref{kernel} which results in
\begin{equation}
\begin{split}
f_{1,k}(r;a_1,\ldots,a_k)=&\frac{(n-k)!}{n!}\partial_r\int_0^\infty \frac{dt}{(1+t)^{n+1}}\\
&\hspace*{-2.5cm}\times\det \left[\begin{array}{c|c} -K(a_e,-rt) & -K(a_e,a_b)  \\\hline
\displaystyle \overset{\ }{\int_0^r} dv\ K(v,-rt)\left(1-\frac{v}{r}\right)^{n-1}  &
   \displaystyle \overset{\ }{\int_0^r} dv\ K(v,a_b)\left(1-\frac{v}{r}\right)^{n-1}-\left(1-\frac{a_b}{r}\right)^{n-1}\Theta(r-a_b) 
  \end{array}\right]
\end{split}
\end{equation}
with the column index $b=1,\ldots,k$ and the row index $e=1,\ldots,k$. After pulling out the minus sign from the first $k$ rows and then interchanging those with the last row, the signs are cancelling. Moreover, we can pull in the integral over $t$ into the first column yielding the first result of \autoref{theo:poly ensemble}. 

When carrying out the derivative in $r$ we arrive at the second expression. To achieve this, one can do a Laplace expansion along the first column (or first row) in the first expression \eqref{eq:1,kpt poly} and integrate by parts in $t$, using the fact that
\begin{equation}
\partial_r K(x,-rt)=\frac{t}{r} \partial_t K(x,-rt).
\end{equation}
The differentiability is guaranteed and the boundary terms vanish because of $n>2$ and the factor of $t$ from the latter equation.

    
 \section{Application to P\'olya Ensembles}\label{sec:app}

To make the expressions more readable so that we gain more insight, we will adopt the following notation for the Mellin transform throughout this section,
\Eq{Mel.def}{
\tilde{f}(s):=\mathcal{M}f(s)=\int_0^\infty u^{s-1}f(u) du
}
and its incomplete Mellin transform
\Eq{incMel.def}{
\tilde{f}_x(s):=\int_0^x u^{s-1}f(u) du, \quad x\geq 0.
}
For instance,  for the Laguerre and Jacobi ensemble the incomplete Mellin transform $\Tilde{w}_x$ take the form
\Eq{incgam.def}{
 \tilde{w}_{{\rm Lag},x}(c+1)=\int_0^x  w_{\rm Lag}(u) u^cdu=\int_0^x   u^{c+\alpha}e^{-u}du
 =\gamma(c+\alpha+1,x)
 }
 and
\begin{equation}\label{incbet.def}
\begin{split}
\tilde{w}_{{\rm Jac},x}(c+1)=&\int_0^x  w_{\rm Jac}(u) u^c du=\int_0^x u^{c+\alpha}(1-u)^{\beta+n-1} du\\
=&\,\mathrm{B}(x,c+\alpha+1,\beta+n),
\end{split}
\end{equation}
respectively. We have employed the lower incomplete Gamma function $\gamma$ and the incomplete Beta function $\mathrm{B}$. 
The weight function for the Laguerre ensemble is given by
\Eq{Lag.def}{
w_{\rm Lag}(x):=x^\alpha e^{-x},\qquad \alpha>-1}
while it is
\Eq{Jac.def}{
w_{\rm Jac}(x):=x^\alpha (1-x)^{\beta+n-1} \Theta(1-x),\qquad  \alpha,\beta>-1} 
for the Jacobi ensemble.
If one wants to have an existing function $q_n$ for this ensemble, cf., \autoref{coro:poly ensemble explicit}, one has to assume $\beta>0$ otherwise the $n$-th derivative of the weight function is not integrable at $x=1$.

\subsection{Proof of \autoref{coro:poly ensemble explicit}}\label{comprho}

To prove this proposition we need the following lemma.

\lem{prop:ker1}{
Let $K$ be the kernel~\eqref{kernel} of a polynomial ensemble. Then,
\Eq{eq:ker1}{
\int_0^\infty  K(x,y) x^k dx=y^k\qquad{\rm for\ all}\ k=0,\ldots,n-1.
}
}

\begin{proof}[Proof of \autoref{prop:ker1}]
Let, $\left\{p_b\right\}_{b=0}^{n-1}$, with $p_b$ a polynomial of degree $b$, and  $\left\{q_b\right\}_{b=0}^{n-1}$ be the bi-orthonormal system given by the polynomial ensemble, i.e.,
\Eq{}{\int_0^\infty dx\ q_b(x)p_c(x)=\delta_{c,b}.}
Any polynomial $f$ of degree less than $n$ has a unique decomposition in the basis $\left\{p_b\right\}_{b=0}^{n-1}$, namely
\Eq{}{
f(x)=\sum_{b=0}^{n-1} \left(\int_0^\infty dt\, q_b(t)f(t)\right) p_b(x).
}
Taking $f(y)=y^k$ for a $k\in \llbracket 0,n-1 \rrbracket$ the result follows.
\end{proof}

\begin{proof}[Proof of \autoref{coro:poly ensemble explicit}]
Starting from the second line of Eq.~\eqref{eq:1,kpt poly} and subtracting $f_{1,0}(r)f_{0,k}(a_1,\ldots,a_k)$ we find
\begin{equation}
\begin{split}
\mathrm{cov}_{1,k}(r;a_1,\ldots,a_k)=&\frac{(n-k)!}{(n-1)!}\int_{0}^{\infty} dt  \det\left[\begin{array}{c|c} 
    	 \displaystyle 0 & \Omega(r,a_c,t) \\ \hline
     K\left(a_b,-rt\right)\quad & K(a_b,a_c) 
\end{array}\right]_{b,c=1}^k\\
=&-\frac{(n-k)!}{(n-1)!}\sum_{l,m=1}^k(-1)^{l+m}\det[K(a_b,a_c)]_{\substack{b\in\llbracket 1,k\rrbracket\setminus\{l\}\\c\in\llbracket 1,k\rrbracket\setminus\{m\}}}\\
&\times\underbrace{\int_{0}^{\infty} dt\ K\left(a_l,-rt\right)\Omega(r,a_m,t)}_{=:-\hat{C}(r;a_l,a_m)}.
\end{split}
\end{equation}
We can rewrite the sum by introducing an auxiliary parameter $\mu$ and understand the expression as a derivative of a determinant, using Jacobi's formula, i.e.,
\begin{equation}
\begin{split}
\mathrm{cov}_{1,k}(r;a_1,\ldots,a_k)=&\frac{(n-k)!}{(n-1)!}\partial_\mu\det\left[\begin{array}{c} K(a_b,a_c) +\mu\,\hat{C}(r;a_b,a_c)
\end{array}\right]_{b,c=1}^k\biggl|_{\mu=0}.
\end{split}
\end{equation}
The expression~\eqref{cov.prop2} is left to show. After performing the change of variables $-t=\tau/(\tau-1)$ and employing the definitions~\eqref{def:phi} and~\eqref{def:omega}, it is
\Eq{eq:cov2}{
\hat{C}(r;a_1,a_2)=\int_0^\infty dv\ T(r,v,a_1)K(v,a_2)-T(r,a_1,a_2),
}
with
\Eq{T.def}{
T(r,v,a) &:=\frac{1}{n^2}\int _{0}^{1}d\tau \frac{1}{v} h\left(\frac{r}{v},\tau\right) K\left(a,\frac{r\tau}{\tau-1}\right),
}
where 
\Eq{}{
 h(x,\tau) =-n \left[ \Psi_0(x^{-1})(1-\tau)^{n-1}+(n+1)\Psi_1(x^{-1})\tau (1-\tau)^{n-1}\right]\Theta(1-x),
}
with $\Psi_0$ and $\Psi_1$ defined in~\eqref{psi.def}. Using the integral representation of the kernel~\eqref{eq:kernel polya} in~\eqref{T.def}, we obtain
\Eq{}{
T(r,v,a) =\frac{1}{n}\int _{0}^{1} \frac{d\tau}{v} h\left(\frac{r}{v},\tau\right) \int_0^1 dt q_n(a t)p_{n-1}\left(\frac{r t\tau}{\tau-1}\right).
}
The absolute integrability in $t$ is given by $t\mapsto q_n(a t)$ while the one in $\tau$ follows from $(1-\tau)^{n-1}|p_{n-1}(rt\tau/(\tau-1))|<\infty$ for all $\tau\in[0,1]$. Thus, we can interchange the two integrals and find
\Eq{}{
 T(r,v,a) =&-\frac{1}{v}\Theta(v-r) \bigg[\Psi_0\left(\frac{v}{r}\right)\int_0^1 dt q_n(a t)\int _{0}^{1}d\tau p_{n-1}\left(\frac{r t\tau}{\tau-1}\right)(1-\tau)^{n-1}\\
& +(n+1)\Psi_1\left(\frac{v}{r}\right) \int_0^1 dt q_n(a t)  \int _{0}^{1}d\tau  p_{n-1}\left(\frac{r t\tau}{\tau-1}\right) \tau (1-\tau)^{n-1}\bigg].\\
}
To shorten the notation, we define
\Eq{P.def}{
P_\gamma(x) := (n+1)^\gamma\int _{0}^{1}d\tau  p_{n-1}\left(\frac{ x \tau}{\tau-1}\right) \tau^\gamma (1-\tau)^{n-1},
\ 
H_\gamma(x,y) :=\int_0^1 dt P_\gamma(x t) q_n(y t)
}
for $\gamma=0,1$. In terms of these functions, we have
\Eq{}{
T(r,v,a) =-\frac{1}{v}\Theta(v-r) \sum_{\gamma=0,1}\Psi_\gamma\left(\frac{v}{r}\right)H_\gamma(r,a)
}
and
\Eq{T.int}{
\int_0^\infty T(r,v,a)K(v,a) dv = -n\int_1^\infty \frac{du}{u}\int_0^1 dt \sum_{\gamma=0,1}\Psi_\gamma(u)H_\gamma(r,a)  q_n(r u t)p_{n-1}(at),}
where we substituted $v=ru$.

To interchange the integrals in~\eqref{T.int}, we need to be more careful as the integration is not absolutely convergent when combined in $u$ and $t$ as the weight function $q_n$, which  usually guarantees the convergence, depends on the product $ut$ which is troublesome when $t$ becomes small while $u$ is large. The idea is to go back to expression~\eqref{kernel} of the kernel in terms of a sum and use the bi-orthonormality relations between  $p_k$ and $q_j$. We note that the function $u \mapsto u^{-1}\Psi_\gamma(u)$ is a polynomial in $u$ of degree $n-1$ for both $\gamma=0,1$. Thence, it can be decomposed in terms of the polynomials $u\mapsto p_j(xu)$ with an arbitrary $x>0$, i.e.,
\Eq{phi.pol}{u^{-1}\Psi_\gamma(u)=\sum_{k=0}^{n-1} a_k(x) p_k(xu).
}
The auxiliary variable $x$ is important since the kernel comes with $x=r$. We can combine this with
\Eq{}{
\int_1^\infty   \frac{du}{u}\Psi_\gamma(u)
K(xu,y)=\int_0^\infty   \frac{du}{u}\Psi_\gamma(u)
K(xu,y)-\int_0^1   \frac{du}{u}\Psi_\gamma(u)
K(xu,y).
}
Due to \autoref{prop:ker1} we have then
\Eq{}{
\int_0^\infty   \frac{du}{u}\Psi_\gamma(u)
K(xu,y)=&\sum_{k=0}^{n-1} a_k(x) \int_0^\infty du\ p_k(xu)K(xu,y)=\sum_{k=0}^{n-1} a_k(x) \frac{p_k(y)}{x}=\frac{1}{y}\Psi_\gamma\left(\frac{y}{x}\right).
}

The integrals in the remaining integration over $u,t\in[0,1]$ can be interchanged now, leading to
\Eq{}{
\int_0^1   \frac{du}{u}\Psi_\gamma(u)
K(xu,y)= n\int_0^1 dt  p_{n-1}(yt) \int_0^1\frac{du}{u}\Psi_\gamma(u)
  q_n(x u t).
}
We combine these considerations with
\Eq{}{
\int_0^\infty \frac{du}{u}\Psi_\gamma(u)q_n(x u t)=0
}
and plug them into~\eqref{T.int} to arrive at
\begin{equation}
\begin{split}
\int_0^\infty  T(r,v,a)K(v,a) dv =&-\sum_{\gamma=0,1}\frac{1}{a}\Psi_\gamma\left(\frac{a}{r}\right)H_\gamma(r,a)\\
&-n\int_0^1dt\sum_{\gamma=0,1}H_\gamma(r,a)p_{n-1}(at)\int_1^\infty\frac{du}{u}\Psi_\gamma(u)q_n(rut).
\end{split}
\end{equation}
We define for $\gamma=0,1$ the functions
\Eq{Q-W.def}{
Q_\gamma(x) := -n\int_1^\infty \frac{du}{u} \Psi_\gamma(u)q_n(x u), 
\qquad
V_\gamma(x,y) := \int_0^1 dt Q_\gamma(x t) p_{n-1}(y t),
}
so that the cross-covariance density~\eqref{eq: cov def} takes the desired form~\eqref{cov.prop}.

What is left is to compute $Q_\gamma$ and $P_\gamma$ using the explicit expressions~\eqref{eq:polya ortho functions} for $p_{n-1}$ and $q_n$.  For the function $P_\gamma$ we start from its definition~\eqref{P.def} and compute
\begin{equation}
\begin{split}
P_\gamma(x) 
=& (n+1)^\gamma\sum_{c=0}^{n-1} \binom{n-1}{c} \frac{ x^c}{\Tilde{w}(c+1)}\int _{0}^{1}d\tau \tau^{\gamma+c} (1-\tau)^{n-1-c}.
\end{split}
\end{equation}
Using the formula for the Beta function, we get
\Eq{}{
P_\gamma(x) =\frac{1}{n} \sum_{c=0}^{n-1}  \frac{ x^c}{\Tilde{w}(c+1)}\left(c+1\right)^\gamma
}
which is, however, only valid for $\gamma=0,1$.
The resulting sum can be identified with the 1-point function of the squared eigenradii for P\'olya ensembles, see~\cite[Lemma 4.1]{Kieburg2016}. Thence, it is
\Eq{P.exp}{
  P_0(x) = \displaystyle \frac{\rho_{\rm EV}(x)}{w(x)}\qquad{\rm and}\qquad P_1(x) = \displaystyle \partial_x \left[ x\frac{\rho_{\rm EV}(x)}{w(x)}\right].
}

For the function $Q_\gamma$, we use the derivative formula~\eqref{eq:polya ortho functions} for $q_n$ and plug it into the definition~\eqref{Q-W.def}. For this aim, we consider the integral
\Eq{}{
\widehat{Q}_\gamma(x) &:=\int_1^\infty du\, q_n(x u) (u-1)^{n-2+\gamma}.
}
Changing the integration variable $v=xu$, we need to compute
\Eq{}{
\widehat{Q}_\gamma(x) =\frac{x^{-(n-1+\gamma)}}{n!}\int_x^\infty dv\, \partial_v^n[v^n w(v)] (v-x)^{n-2+\gamma}.
}
We can now integrate by parts $n-2+\gamma$ times. As the boundary terms are vanishing, this yields
\Eq{}{
\widehat{Q}_0(x) &= \frac{(-1)^{n-1}}{n(n-1)} [n w(x) +x \partial_x w(x)],\qquad \widehat{Q}_1(x) = \displaystyle \frac{(-1)^{n} }{n}  w(x).
}
The functions $Q_\gamma$ expressed in terms of $\widehat{Q}_\gamma$ are
\begin{equation}\label{Q.exp}
\begin{split}
Q_0(x) =& (-1)^{n-1}n [n \widehat{Q}_1(x)+ (n-1) \widehat{Q}_0(x)]= x \partial_x w(x) ,\\
Q_1(x) =&(-1)^{n-2}n  \widehat{Q}_1(x)= w(x) .
\end{split}
\end{equation}
Plugging this into~\eqref{Q-W.def} concludes the proof. 
\end{proof}

Almost all the integrals involved in the formulation of \autoref{theo:poly ensemble} can be carried out, yielding an expression of the cross-covariance density which is very efficient for numerical evaluations. The only remaining integration is hidden in the incomplete Mellin transform of $w$. We will use \autoref{coro:poly ensemble explicit} to prove the following corollary for $\mathrm{cov}(r;a)=\hat{C}(r;a,a)$.

\coro{lem:poly ensemble explicit}{
Let $n\in \mathbb{N}$, $n>2$.  With the same assumptions and notations as in \autoref{theo:poly ensemble} and choosing a P\'olya ensemble associated to the weight function $w$, the cross-covariance density can be cast into the form 
\begin{equation}\label{eq:poly ensemble explicit}
\begin{split}
\mathrm{cov}(r;a)=&\frac{\Theta(r-a)}{nar} \sum_{c=0}^{n-1}  \frac{ (r/a)^c}{\Tilde{w}(c+1)}  \Tilde{q}_{n,a}(c+1)\left(1-\frac{a}{r}\right)^{n-2}\left[    (n-c-1)\frac{a}{r}+c\right]\\
&\hspace*{-1.5cm}-\frac{1}{na} \sum_{c=0}^{n-1}  \frac{ (r/a)^c}{\Tilde{w}(c+1)} \Tilde{q}_{n,a}(c+1)\left[w(r) p_{n-1}(a)+ \sum_{j=0}^{n-1} \binom{n-1}{j}(c-j)  \frac{\Tilde{w}_r(j+1)}{\Tilde{w}(j+1)} \frac{(-a)^j}{r^{j+1}} \right].
\end{split}
\end{equation}
We recall the definitions and notations of the Mellin transform~\eqref{Mel.def} and the incomplete Mellin transform~\eqref{incMel.def}.
}

\begin{proof}[Proof of \autoref{lem:poly ensemble explicit}]

Considering the proof of \autoref{coro:poly ensemble explicit}, we still need to compute $V_\gamma$ and $H_\gamma$. Using the sum expression~\eqref{biorthogonal} of $p_{n-1}$ and the explicit expression~\eqref{Q.exp} of $Q_\gamma$ with the change of variables $u=xt$, one gets
\Eq{}{
V_1(x,y)&=\sum_{j=0}^{n-1} \binom{n-1}{j}  \frac{\left(-y/x\right)^j}{\Tilde{w}(j+1)} \frac{\Tilde{w}_x(j+1)}{x},
}
while for $V_0$ it is
\Eq{eq:V_0}{
V_0(x,y) &= \sum_{j=0}^{n-1} \binom{n-1}{j}  \frac{\left(-y/x\right)^j}{\Tilde{w}(j+1)} \frac{1}{x}\int_0^x  du\  u^{j+1}\partial_u w(u).
}
One can proceed with integration by parts to arrive at
\Eq{}{
V_0(x,y) &= w(x) p_{n-1}(y)-\partial_y\left[y V_1(x,y)\right].
}

Similarly to $V_\gamma$, one can also compute $H_\gamma$ by exploiting~\eqref{P.exp} for $P_\gamma$ and the derivative expression~\eqref{biorthogonal} of $q_n$. Then, we obtain for $\gamma=0,1$ 
\Eq{}{
H_\gamma(x,y) &=  \frac{1}{n!} \int_0^1 dt  \sum_{c=0}^{n-1}  \frac{ (xt)^c}{\Tilde{w}(c+1)}\left(c+1\right)^\gamma \partial_{(yt)}^n[(yt)^n w(yt)].
}
Changing to $u=yt$, it is
\Eq{}{
H_\gamma(x,y) &=  \frac{1}{n!} \sum_{c=0}^{n-1}  \frac{ (x/y)^c}{\Tilde{w}(c+1)} \left(c+1\right)^\gamma \frac{1}{y}\overbrace{\int_0^y du\  u^c \partial_u^n[u^n w(u)]}^{=n!\,\Tilde{q}_{n,y}(c+1)}.
}
When plugging these explicit expressions for $H_\gamma$ and $V_\gamma$ into~\eqref{cov.prop} one gets the claim of \autoref{lem:poly ensemble explicit}. 
\end{proof}

We would like to point out that the remaining integral is nothing else than the incomplete Mellin transform of $q_n$ which we can express in terms of a sum after integrating by parts and collecting all boundary terms, 
\Eq{}{
n!\,\Tilde{q}_{n,y}(c+1)=\int_0^y  u^c \partial_u^n[u^n w(u)] du=\sum_{p=0}^c \binom{c}{p}p!(-1)^p y^{c-p} \partial_y^{n-1-p}[y^n w(y)].
}
\begin{remark}
Equation~\eqref{eq:poly ensemble explicit} can also be cast into the factorized form
\Eq{}{
\mathrm{cov}(x;y)=&\frac{1}{n xy}\sum_{c,j=0}^{n-1} \binom{n-1}{j} (-1)^j \left(\frac{y}{x}\right)^{j-c} \frac{\Tilde{q}_{n,y}(c+1)}{\Tilde{w}(c+1)}\Bigg[(c+1)\left(\Theta(x-y)-\frac{\Tilde{w}_x(j+1)}{\Tilde{w}(j+1)} \right)\\
&- \Theta(x-y) \left(\frac{x-ny}{x-y}\right)+\frac{(j+1)\tilde{w}_x(j+1)-x^{j+1}w(x)}{\Tilde{w}(j+1)}\Bigg], 
}
where we used the Newton binomial formula for the term $(1-y/x)^{n-2}$ and plugging the identity
\Eq{}{
\int_0^x  u^{j+1}\partial_u w(u)  du=x^{j+1}w(x)-(j+1)\tilde{w}_x(j+1)
}
in~\eqref{eq:V_0}.
\end{remark}

\subsection{Proof of \autoref{coro:deriv Polya}}\label{proof coro deriv polya}

For $n=2$, the explicit formula in \autoref{prop:rho n=2}, clearly shows that the $1,1$- and $1,2$-point correlation functions are discontinuous along the lines $r=a_1$ and $r=a_2$ which is reflected in the Heaviside step functions.

For $n>2$, we consider the expression of the $1,k$-point function given in \autoref{theo:poly ensemble}. It is clear that $\rho_{\rm EV}$ ($(1,1)$ entry of the determinant) and the kernel $K$ are continuous and smooth on the support $\sigma$ for a P\'olya ensemble because of their dependence on the  smooth weight function $w \in \mathrm{C}^{\infty}(\sigma)$, see~\eqref{eq:kernel polya} and~\eqref{P.exp}. The smoothness of the kernel $K$ in its first argument is due to the compactness of the second integral expression in~\eqref{eq:kernel polya} and integrability of $q_n$ on the support  $\sigma$ so that the Leibniz integral rule applies. 

In summary, the reduced differentiability of $f_{1,k}$ must be inherited from the function
\Eq{eq:fct diff}{
(r,v,a)\mapsto \frac{\Theta(r-v)}{v}\int_0^\infty dt\, \varphi\left(\frac{v}{r},t\right) K(a,-rt)
}
see~\eqref{T.def}. Since $K$ is a polynomial of degree $n-1$ in its second entry the integral over $t$ gives a smooth function on the support $\sigma$ in all three variables $r$, $v$ and $t$. Furthermore, this integral does not vanish at almost all $v=r=a$ because of
\begin{equation}
\begin{split}
&\lim_{v,r\to a}\int_0^\infty dt\, (1+t)^{-(n+2)}\left[\left(1-\frac{v}{n r}\right)(1+t)-\left(1-\frac{v}{r}\right)\left(1+\tfrac{1}{n}\right) \right] K(a,-rt)\\
=&\frac{n-1}{n} \int_0^\infty dt\, (1+t)^{-(n+1)} \int_0^1 ds\ q_n(as)p_{n-1}\left(-ast\right)\\
=&\frac{n-1}{n}\int_0^1 ds\ q_n(as)P_0(as)\\
=&\frac{n-1}{an}\int_0^a ds\ q_n(s)\frac{\rho_{\rm EV}(s)}{w(s)}.
\end{split}
\end{equation}
In the first equality we have exploited the integral representation~\eqref{eq:kernel polya} of the kernel, and in the second equality we could interchange the two integrals as they are absolutely integrable due to the polynomial nature of $p_{n-1}$. Additionally, we have employed Eqs.~\eqref{P.def} and~\eqref{P.exp}.
As can be readily checked the derivative in $a$ of $a$ times this integral is equal to $(n-1)q_n(a)\rho_{\rm EV}(a)/w(a)$ which is evidently non-linear implying a non-constant behaviour of this integral when $a\in\sigma$.

In conclusion of this discussion, the function \eqref{eq:fct diff} is $(n-3)$-times continuous differentiable at $r=a$ and its $(n-2)$nd derivative is discontinuous along this line because of the factor $\Theta(a-r)\left(a/r-1\right)^{n-2}$. 

The question is whether the integral
\begin{equation}
\begin{split}
\int_0^r \frac{dv}{v}\int_0^\infty dt\, \varphi\left(\frac{v}{r},t\right) K(a_1,-rt) K(v,a_2)
\end{split}
\end{equation}
may change this differentiability. The point is that the integrand is smooth in $(r,a_1,a_2)\in\sigma^3$ and $(v,t)\in(0,r)\times\mathbb{R}_+$ and absolutely integrable in $v$ and $t$. Thus, the Leibniz integral rule tells us that after integrating $v$ and $t$ it remains smooth in $r$ and $a$. In summary, this tells us that the cross-covariance density~\eqref{eq: cov def} is a sum of a smooth function and an $(n-3)$-times continuous differentiable function on $\sigma^{k+1}$ which has been the statement of the corollary.

\subsection{Examples: Laguerre and Jacobi Ensembles}\label{Example: Orthogonal ensembles}

There are classical random matrix ensembles which are P\'olya ensembles. For instance, the Laguerre ensemble (also known as Ginibre ensemble) is one, see~\cite[Examples 3.4]{Kieburg2016}. The weight function $w=w_{\rm Lag}$ is given by~\eqref{Lag.def} which corresponds to the joint probability function of the squared singular values
\begin{equation}
f_{\rm SV}(a)=\frac{1}{n!}\left(\prod_{j=0}^{n-1}\frac{1}{j!\,\Gamma(\alpha+j+1)}\right)\Delta_n^2(a)\prod_{j=1}^n a_j^\alpha e^{-a_j}.
\end{equation}
Employing the Laguerre polynomials
\Eq{Lag.pol.def}{
 L_{j}^{(\alpha )}(x):=\sum _{k=0}^j\binom{j+\alpha}{j-k} \frac {(-x)^{k}}{k!},
}
we have explicit expressions for the bi-orthonormal set of functions composing the kernel of the determinantal point process,
\begin{equation}
\begin{split}
p_{j}(x)=\frac{j!}{\Gamma(j+\alpha+1)} L_{j}^{(\alpha )}(x)\qquad{\rm and}\qquad q_j(x)= L_{j}^{(\alpha )}(x)x^\alpha e^{-x}.
\end{split}
\end{equation}
It is well-known~\cite[Corollary 5.2]{Kuijlaars2014a} that the kernel can be expressed in terms of the one-fold integral
\begin{equation}
K(x,y)=\frac{n!}{\Gamma(n+\alpha)}\int_0^1dt\ L_{n-1}^{(\alpha )}(yt)L_{n}^{(\alpha )}(xt)(xt)^\alpha e^{-xt}.
\end{equation}
What is new is the cross-covariance density of a squared singular value $a$ with a squared eigenradius $r$. We make use of \autoref{lem:poly ensemble explicit} and plug in the Mellin transform and incomplete Mellin transform of  $w=w_{\rm Lag}$,
\begin{equation}
\begin{split}
\tilde{w}_{\rm Lag}(c+1)=\Gamma(c+\alpha+1)\qquad{\rm and}\qquad \tilde{w}_{{\rm Lag},x}(c+1)=\gamma(c+\alpha+1,x),
\end{split}
\end{equation}
see~\eqref{incgam.def}, as well as the incomplete Mellin transform of the weight $q_n$ which is essentially a hypergeometric function
\begin{equation}
\begin{split}
\tilde{q}_{n,x}(c+1)=&\frac{1}{n!}\int_0^x du\, u^c \partial_u^n[u^{n+\alpha}e^{-u}]\\
=&\frac{x^{\alpha+c+1}}{n!}\sum_{j=0}^\infty\frac{\Gamma(n+\alpha+j+1)}{\Gamma(\alpha+j+1)(\alpha+c+j+1)}\frac{(-x)^j}{j!}\\
=&\frac{\Gamma(n+\alpha+1)}{n!\,(\alpha+c+1)\Gamma(\alpha+1)} x^{\alpha+c+1}\\
&\times { _2F_2}(\alpha+c+1,n+\alpha+1;\alpha+1,\alpha+c+2;-x).
\end{split}
\end{equation}
Unfortunately, we were unable to simplify the expression further for this case.

Another classical ensemble falling into the class of P\'olya ensembles is the Jacobi ensemble (also known as truncated unitary ensemble). Using the weight function $w=w_{\rm Jac}$ in~\eqref{Jac.def}, one can find the joint probability function of the squared singular values~\cite[Examples 3.4]{Kieburg2016}
\begin{equation}
f_{\rm SV}(a)=\frac{1}{n!}\left(\prod_{j=0}^{n-1}\frac{\Gamma(n+\alpha+\beta+j+1)}{j!\,\Gamma(\alpha+j+1)\Gamma(n+\beta)}\right)\Delta_n^2(a)\prod_{j=1}^n a_j^\alpha (1-a_j)^\beta\Theta(1-a_j).
\end{equation}
This time, the Jacobi polynomials are involved,
\begin{equation}
\begin{split}
P_{j}^{(\alpha ,\beta )}(x):=&\frac{1}{2^j j!\,(1-x)^\alpha(1+x)^\beta}(-\partial_x)^j[(1-x)^{j+\alpha}(1+x)^{j+\beta}]\\
=& \frac{\Gamma(j+\alpha+1)}{j!\,\Gamma(j+\alpha+\beta+1)}\sum _{k=0}^{j}\binom{j}{k}\frac{\Gamma(j+\alpha+\beta+k+1)}{\Gamma(\alpha+k+1)}\left({\frac {x-1}{2}}\right)^{k}.
\end{split}
\end{equation}
It can be shown that
\begin{equation}
\begin{split}
p_{j}(x)=&\frac{j!\,\Gamma(n+\alpha+\beta+1)}{\Gamma(j+\alpha+1)\Gamma(n+\beta-j)} P_{j}^{(\alpha,\beta+n-j)}(1-2x),\\
 q_j(x)=&P_{j}^{(\alpha,\beta+n-j-1 )}(1-2x)\,x^\alpha (1-x)^{\beta+n-j-1},
\end{split}
\end{equation}
because of the Mellin transform of $w=w_{\rm Jac}$
\begin{equation}
\begin{split}
\tilde{w}_{\rm Jac}(c+1)=\frac{\Gamma(c+\alpha+1)\Gamma(n+\beta)}{\Gamma(n+\alpha+\beta+c+1)}.
\end{split}
\end{equation}
We note that the resulting Jacobi polynomials are not the standard ones when approaching this ensemble with orthogonal polynomials. The reason is that we constructed those via bi-orthonormality which is certainly also allowed. The kernel, thus, takes the form
\begin{equation}
K(x,y)=\frac{n!\,\Gamma(n+\alpha+\beta+1)}{\Gamma(n+\alpha)\Gamma(\beta-1)}\int_0^1dt\ P_{n-1}^{(\alpha,\beta+1)}(1-2yt)P_{n}^{(\alpha,\beta-1 )}(1-2xt)\,(xt)^\alpha (1-xt)^{\beta-1}.
\end{equation}
see also~\cite[Proposition 2.7]{Kieburg2015} for $r=1$.

The incomplete Mellin transforms needed for \autoref{lem:poly ensemble explicit} are the incomplete Beta function~\eqref{incbet.def} and, anew, a hypergeometric function
\begin{equation}
\begin{split}
\tilde{q}_{n,x}(c+1)=& \frac{1}{n!}\int_0^x du\, u^c \partial_u^n[u^{n+\alpha}(1-u)^{\beta+n-1}]\\
=&\frac{x^{\alpha+c+1}}{n!}\sum_{j=0}^\infty\binom{n+\beta-1}{j}\frac{\Gamma(n+\alpha+j+1)}{\Gamma(\alpha+j+1)(\alpha+c+j+1)}(-x)^j\\
=&\frac{\Gamma(n+\alpha+1)}{n!\,\Gamma(\alpha+1)(\alpha+c+1)} x^{\alpha+c+1}\\
&\times { _3F_2}(\alpha+c+1,n+\alpha+1,1-n-\beta;\alpha+1,\alpha+c+2;x).
\end{split}
\end{equation}
In this calculation we have used $x\in(0,1)$ as this is the support of the squared eigenradii and squared singular values.

\section{Discussion}\label{Discussion}

We studied cross-correlation functions between the singular values and eigenvalues of arbitrary bi-invariant complex square matrices and found results for the $1,k$-point correlation function $f_{1,k}$, see \autoref{theo:1,kpt}, and for  the $1$-point function $\rho_{\rm EV}$ of the eigenradii, see \autoref{theo:1pt}. These formulas  drastically simplify in the case the bi-unitarily invariant ensemble is a polynomial ensemble, see \autoref{rem:bui ensemble}, \autoref{theo:poly ensemble} and \autoref{rem:cov}. The obtained formula enables us to define and identify the $1,k$-cross-covariance density function between one eigenradius and $k$ singular values, which captures the interaction between the two kind of variables. Note that the interaction could also be studied by defining analogues of cluster functions used in the physics literature~\cite[Eq.(5)]{Dyson1962}, involving, in this case, the $j,k$-point correlation functions.

The simplification of the results increases further for P\'olya ensembles~\cite{Kieburg2016}; see \autoref{coro:poly ensemble explicit}. Although the proof for $n=2$ is very different, one can check that \autoref{coro:poly ensemble explicit} agrees with \autoref{prop:rho n=2} for an arbitrary polynomial ensemble, especially Eq.~\eqref{Polyn.n=2}.

Let us underline that the only known formula for $\rho_{\rm EV}$ for bi-unitarily ensembles in terms of quantities of the squared singular values $f_{\rm SV}$ was, so far, only for P\'olya ensembles; see~\cite[Eq.(4.4)]{Kieburg2016}. We found the generalisation~\eqref{eq:1pointreal} to polynomial ensembles which is a much larger class than P\'olya ensembles.

The formulas~\eqref{eq:1,kpt poly} and~\eqref{cov.prop} for the $1,k$-point correlation function were not known even  for P\'olya ensembles. Following \autoref{rem:cov}, they can be expressed in terms of the $1,k$-cross-covariance density.

A generalisation of our results for $f_{1,k}$ to an arbitrary $j,k$-point correlation function $f_{j,k}$ for $j$ squared eigenradii and $k$ squared singular values with $j>1$ is certainly desirable, however, it will be challenging to obtain explicit results due to the Vandermonde determinant $\Delta_n(s)$ in the SEV transform~\eqref{eq:SEV}. Indeed, this term couples all the poles in the computation of the $n$-fold complex integral and in general the regularisation function cannot be dropped. 
Another issue is that there is a trade-off between getting formulas for more eigenradii and having explicit and compact formulas as each fixed eigenradius comes with an additional complex contour integral.

Another goal of our study has been to find all conditional marginal measures of the eigenvalues when fixing the singular values. We have been successful for the level density of the eigenradii conditioned to fixed singular values, see \autoref{cor:cond.dens}. The $n$-point correlation function for the eigenradii conditioned to singular values is more complicated though we have a conjecture. Indeed, $2$-point correlation function for $n=2$ can be read off from \autoref{prop:rho n=2} summarised in the following corollary.

\coro{prop:conditional meas n=2}{
For $n=2$, the conditional probability measure of the squared eigenradii $r_1,r_2$ under the condition of the squared singular values $a_1,a_2$ is given by
\Eq{}{
d\mu_{2,2}(r_1,r_2|a_1,a_2)=&\Theta\left(\max\{a_1,a_2\}-\max\{r_1,r_2\}\right)\Theta\left(\min\{r_1,r_2\}-\min\{a_1,a_2\} \right)\\
&\times \frac{r_1+r_2}{2|a_1-a_2|}\, \delta(r_1 r_2-a_1 a_2) dr_1dr_2.
}
}

We recall that the Dirac distribution $\delta(r_1 r_2-a_1 a_2)$ reflects the identity~\eqref{eq:prodsev} while the Heaviside step functions encode the inequalities~\eqref{eq:Weyl1}. Considering the SEV transform~\eqref{eq:SEV} one can conjecture that the conditional joint probability density of the eigenvalues is given in a distributional way by 
\begin{equation}
\begin{split}
d\mu_{n,n}(z|a)=&\frac{\prod_{j=0}^{n-1} j!}{(n!)^2 \pi^n} |\Delta_n(z)|^2dz\lim_{\varepsilon \to 0} \int_{\mathcal{C}(n)}\left[\prod_{k=1}^n \frac{ds_k}{2\pi i}\zeta(\varepsilon \Im{s_k})\right] \\
 & \times  \mathrm{Perm}[|z_b|^{-2 s_c}  ]_{b,c=1}^n\frac{\det [a_b^{s_c-1}  ]_{b,c=1}^n }{\Delta_n(s)\Delta_n(a)} .
\end{split}
\end{equation}
The condition~\eqref{eq:prodsev} is then still encoded in the integral over $s$ and the limit $\varepsilon\to0$. At the moment, one needs to take this formula with caution as a rigorous mathematical proof is missing for this conjecture because the limit $\varepsilon\to0$ cannot exist point-wise but has to be understood in a distributional way.

Reversing the conditional probability measure, meaning finding the measure for the squared singular by conditioning the eigenvalues,  is expected to be more challenging when considering the inverse SEV transform, see~\cite[Eq.~(3.4)]{Kieburg2016}. 

\begin{figure}[t!]
        \centerline{\includegraphics[width=\textwidth]{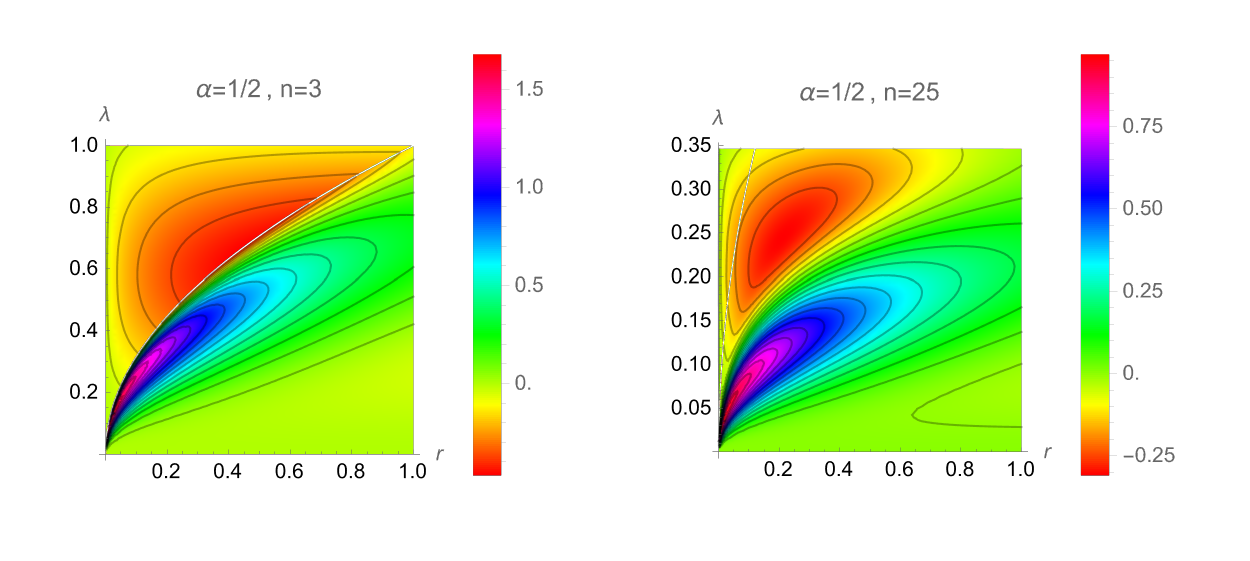}} 
        \centerline{\includegraphics[width=\textwidth]{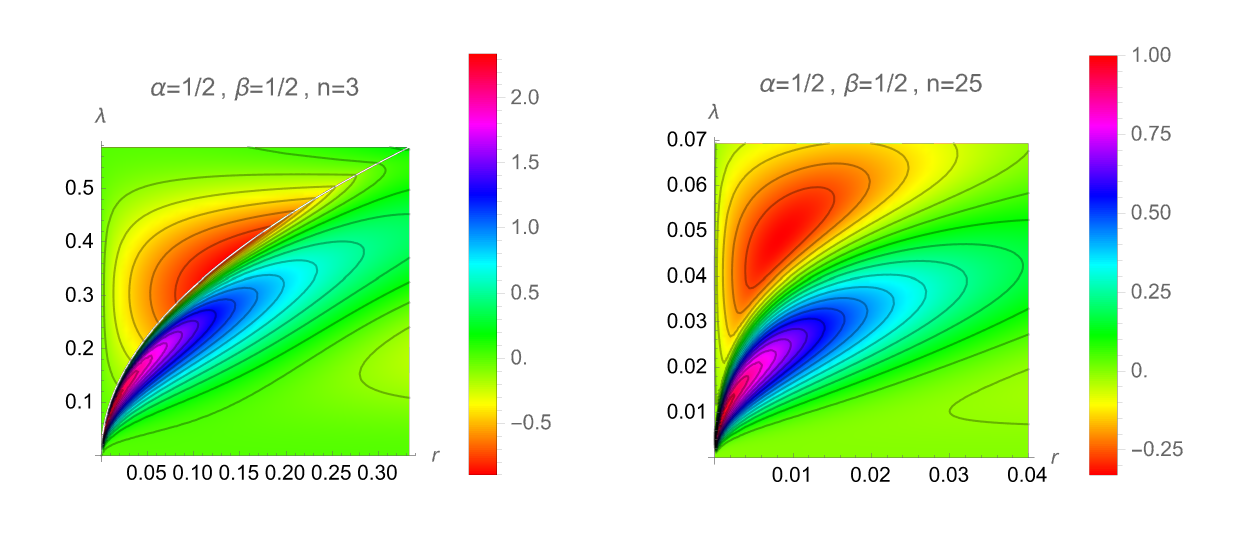}} 
     \caption{\footnotesize Contour plots of the cross-covariance density  $2\lambda \mathrm{cov}(r;\lambda^2)$ for the Laguerre ensemble with $\alpha=1/2$ (upper row) and for the Jacobi ensemble with $\alpha=\beta=1/2$ for the matrix dimensions $n=3$ (left column) and $n=25$ (right column). For the Laguerre ensemble we rescaled by $n^{3/2}$. See Eqs.~\eqref{Lag.def} and~\eqref{Jac.def} for the corresponding weight functions and their parameter dependence in $\alpha$ and $\beta$, respectively. We zoomed into the scale of the local mean level spacing for the singular values $\lambda$ and squared eigenradii $r$. The impact of the line $\lambda=\sqrt{r}$, where the cross-covariance density is not smooth, is still visible for $n=3$ but becomes hard to detect for larger matrix size. } 
    \label{fig:hard}
\end{figure}

Our results for polynomial ensembles  and especially for P\'olya ensembles open a pathway to study the large $n$ limit for cross-correlations of eigenvalues and singular values. For instance, it allows to analyse finite $n$ corrections of  Feinberg-Zee's Single Ring Theorem~\cite{Feinberg1997}, rigorously proven in~\cite{Guionnet2009}, for which some hard-to-check condition has then been lifted by~\cite{Rudelson2014}. The same holds true for the Haagerup-Larsen Theorem~\cite{Haagerup2000} which relates the probability density of the eigenradii with those of the singular values. The finite $n$ analogue is encoded in \autoref{theo:1pt} or more explicitly for polynomial ensembles in~\eqref{eq:1pointreal}. In particular, formula~\eqref{eq:1pointreal} will enable us to understand how the asymptotics described by these theorems are approached. 

Recently, a deformed Single Ring Theorem has been proven in~\cite{Ho2023}. It arises when the bi-unitary invariance is perturbed. Some kind of such perturbations might be also possible to study with  our methods, as already pointed out in~\cite[Sec.~3.3]{Kieburg2016}.

There will also be interesting local spectral statistics of the cross-correlations between the complex eigenvalues and singular value. For instance the condition~\eqref{eq:boundsvev}, implies non-trivial correlations around common edges. In the case of the common hard edge at the origin, it can already be seen in \autoref{fig:hard}, where we show contour plots of the cross-covariance density $2\lambda{\rm cov}(r;\lambda^2)$ for a Laguerre and a Jacobi ensemble with the same parameter $\alpha$, see Eqs.~\eqref{Lag.def} and~\eqref{Jac.def}, for which we know they share the same hard edge statistics of the singular values. We have chosen the singular values $\lambda$ instead of the squared singular values $a=\lambda^2$ as, then, the scaling in $n$ of the mean level spacing agrees with the one of the squared eigenradii. 

The probabilistic interpretation of the cross-covariance density, shown in \autoref{fig:hard}, is as follows. Negative regions translate into some kind of repulsion between the eigenradius and the singular value. It is less likely to find a correlated pair of variables in this region, compared to the case where they would be independent. Similarly, positive regions translate into some kind of attraction between the eigenradius and the singular value. Null regions are simply neutral, the variables are almost uncorrelated there. It is, however, important to stress that, like every covariance, the cross-covariance density is weighted with the likelihood of having at all a squared singular value or a squared eigenradius in the respective intervals. The less likely it is to find one of the variables in a specific interval, the smaller the cross-covariances becomes.

Note that there is no white regions in the plots, despite Weyl's inequalities \eqref{eq:Weyl1}, \eqref{eq:Weyl2} and \eqref{eq:boundsvev}, because the eigenradii and singular values are not ordered here. Indeed, the $j,k$-correlation measures are fully symmetric in the eigenradii and singular values (cf. \autoref{def:j,k pt meas}), i.e. averaged on all possible orderings.

The contour plots for the matrix size $n=25$ show two things: Firstly, non-trivial (non-factorising) spectral statistics seem to emerge in the hard edge limit.  Secondly, these statistics seem to be shared by different ensembles as they appear to match for the two ensembles, apart from a scaling. We investigate this in a follow-up work \cite{Allard2025}.

\section*{Acknowledgement}
We want to thank Arno Kuijlaars for his advices and fruitful discussions. MA acknowledge financial support by the International Research Training Group (IRTG) between the University of Melbourne and KU Leuven and MK has been supported by the Australian Research Council via the Discovery Project grant DP210102887.

\bibliography{main}

\end{document}